\documentclass[12pt]{article}

\title{Diophantine Approximation on varieties I: Algebraic distance and metric B\'ezout Theorem}

\author{Heinrich Massold}

\usepackage{bezier,amsfonts}

\usepackage{amsmath,amsfonts,bbm,amssymb}

\newtheorem{Satz}{}[section]

\newcommand{\satz}[1]{\vspace{2mm} \begin{Satz}{\bf #1}}

\newcommand{\proof}{\vspace{3mm} {\sc Proof}\hspace{2mm}}

\newcommand{\la}{\langle}

\newcommand{\ra}{\rangle}

\newcommand{\di}{{\mbox{div}}}

\newcommand{\R}{{\mathbbm R}}

\newcommand{\Pe}{{\mathbbm P}}

\newcommand{\Z}{\mathbbm{Z}}

\newcommand{\N}{\mathbbm{N}}

\newcommand{\CC}{{\cal C}}

\newcommand{\CQ}{{\cal Q}}

\newcommand{\CG}{{\cal G}}

\newcommand{\CD}{{\cal D}}

\newcommand{\CO}{{\cal O}}

\newcommand{\CX}{{\cal X}}

\newcommand{\CY}{{\cal Y}}

\newcommand{\CZ}{{\cal Z}}

\newcommand{\CW}{{\cal W}}

\newcommand{\CL}{{\cal L}}

\newcommand{\CF}{{\cal F}}

\newcommand{\Q}{{\mathbbm Q}}

\newcommand{\C}{{\mathbbm C}}

\newcommand{\spec}{\mbox{Spec} \:}

\begin{document}

\parindent0mm

\maketitle

\thispagestyle{empty}

\begin{abstract} 
For two properly intersecting effective cycles in projective space $X,Y$, and their intersection product $Z$, the metric B\'ezout Theorem relates the degrees, heights of $X,Y$, and $Z$, as well as their distances and algebraic distances to a given point $\theta$. Applications of this Theorem are in the area of Diophantine Approximation, giving estimates for approximation properties of $Z$ with respect to $\theta$ against the ones of $X$, and $Y$.\end{abstract}

\tableofcontents

\section{Introduction}

%This is the first part in a series of papers on approximations of
%nonalgebraic points on quasi projective varieties by algebraic points
%and the deduction of algebraic independence criteria from these.

Let $k$ be a number field with ring of integers $\CO_k$ and
$\CX$ a flat, integral,
quasi projective scheme of finite type over $\CO_k$. The dimension of 
$\CX$ will be denoted by $t+1$, and the base extension of $\CX$ to $k$ by $X$.
Let further $\bar{\CL}$ be an ample metrized line bundle on some projective
closure $\bar{\CX}$ of $\CX$, and $\sigma: k \hookrightarrow ´\C$ a fixed
embedding which gives rise to a base extension
$\psi:X_\sigma = \CX \times_{\spec \; \CO_k} \spec \; \C_\sigma \to \CX$.
On $X_\sigma(\C)$ choose any metric $|\cdot,\cdot|$ that induces the
usual topology on $X_\sigma(\C)$.

For a point $\alpha \in \CX(\bar{k})$ denote 
$\kappa(\alpha)=\CO_{\CX,\alpha}/m_{\alpha}$
the residue field of $\alpha$ which is a finite extension of $k$.
The stalk $\CL_{\alpha}$ of $\CL$ at $\alpha$ is a one dimensional
$\CO_{\CX,\alpha}$-module, thus $\CL_{\alpha}/m_\alpha \CL_\alpha$ is a one dimensional
$\kappa(\alpha)$- vector space. Hence, with $f_0 \in \Gamma(\bar{\CX},\CL)$ 
such that $f_0|_\alpha \notin m_\alpha \CL_\alpha$, one gets a $k$-linear map
\[ r_D : \Gamma(\CX,\CL^{\otimes D}) \to \kappa(\alpha), \quad
   f \mapsto a \quad \mbox{such that} \quad 
   \bar{f}_\alpha  a (\bar{f}_0)_\alpha^{\otimes D} \in m_\alpha, \] 
whose image does not depend on the choice of $f_0$.

Define $X_D(\bar{k})$ as the set of $\alpha \in \CX(\bar{k})$ such that
$r_D$ is a surjection onto $\kappa(\alpha)$, and $X_{D,H}(\bar{k})$ as the
set of $\alpha \in \CX(\bar{k})$ such that the image of the global sections
$f \in \Gamma(\CX,\CL^{\otimes D})$ with $\log |f| =\log |f|_{L^2(X)} \leq H$ under
$r_D$ spans $\kappa(\alpha)$ as a $k$-vector space.

To a point $\theta \in X_\sigma(\C)$ we attach several numbers: 
With $x=\psi(\theta)$, set $\kappa(x) = \CO_{\CX,x}/m_x$ which is
a finitely generated extension of $k$, and set
$t(\theta)\in \N$ equal to the transcendence degree  of $\kappa(x)$
which equals the dimension of the closure of $x$.

The other numbers attached to $\theta$ are defined via the approximability
of $\theta$ by algebraic points and hypersurfaces.
For $D \in \N$ a $D$-frame $F$ is a pair $(\alpha,f)$ consisting of a
point $\alpha \in X_D(\bar{k})$ together with a nonzero 
$f \in \Gamma(\bar{\CX},\CL^{\otimes D})$ such that $f_\alpha \neq 0$. 
For $D \in \N,H \in \R^{\geq 0}$ a $(D,H)$-frame
$F$ consists of an $\alpha \in X_{D,H}(\bar{D})$ together with a nonzero
$f \in \Gamma(\CX,\CL^{\otimes D})$ of logarithmic length at most $H$
such that $f_\alpha \neq 0$.
Denote by $\CF_D, \CF_{D,H}$ the set of $D$-frames and 
$(D,H)$-frames respectively. 

For $F=(\alpha,f)$ an element of $\CF_D$ or $\CF_{D,H}$ and 
$\theta \in X(\C_\sigma)$ a nonalgebraic point, define
\[ D(F,\theta) := \max(\log |\alpha,\theta|,\log |f_\theta|), \quad
   h(\alpha) \quad \mbox{the hight of $\alpha$}, \]
and further the lower and upper approximational degrees
\[ \underline{t}_1(\theta) := \sup \left\{ s \in \R | 
    \limsup_{D \to \infty} \left(-\inf_{F=(\alpha,f) \in \CF_D} \left\{
            \frac{D(F,\theta)}{(h(\alpha)+\deg \alpha) (\deg \alpha)^s} 
             \right\}     \right)
              =\infty\right\} \in \R^{\geq 0}, \]
\[ \underline{t}_2(\theta) := \sup \left\{ s \in \R | 
    \limsup_{D \to \infty} \left(-\inf_{F=(\alpha,f) \in \CF_D} \left\{
            \frac{D(F,\theta)}{(D+\log |f|)D^s} \right\}
                   \right)
              =\infty\right\} \in \R^{\geq 0}, \]
\[ \bar{t}_2(\theta) := \sup \left\{ s \in \R | 
    \liminf_{D \to \infty} \left(-\inf_{(F=(\alpha,f) \in \CF_D} \left\{
            \frac{D(F,\theta)}{(D+\log |f|)D^s} \right\}
                   \right)
              =\infty\right\} \in \R^{\geq 0}, \]

and for $a \in \R^+$ the upper and lower approximational $a$-degree
\[ \underline{t}_1^a(\theta) := \sup \left\{ s \in \R | 
    \limsup_{D \to \infty} \left(-\inf_{F=(\alpha,f) \in \CF_{D,aD}} \left\{
            \frac{D(F,\theta)}{a(\deg \alpha)^{s+1}} 
             \right\}     \right)
              =\infty\right\} \in \R^{\geq 0}, \]
\[ \underline{t}_2^a(\theta) := \sup \left\{ s \in \R | 
    \limsup_{D \to \infty} \left(-\inf_{F=(\alpha,f) \in \CF_{D,aD}} \left\{
            \frac{D(F,\theta)}{aD^{s+1}} \right\}
                   \right)
              =\infty\right\} \in \R^{\geq 0}, \]
\[ \bar{t}_2^a(\theta) := \sup \left\{ s \in \R | 
    \liminf_{D \to \infty} \left(-\inf_{(F=(\alpha,f) \in \CF_{D,aD}} \left\{
            \frac{D(F,\theta)}{D^{s+1}} \right\}
                   \right)
              =\infty\right\} \in \R^{\geq 0}, \]

The main objective of this series of papers is to proof the 

\satz{Theorem} \label{approx}
\begin{enumerate}
Let $\CX$ be any flat, integral, quasi projective scheme over $\CO_k$, choose 
an embedding $\sigma:k\to\C$, let $\bar{\CL}$ be a metrized line bundle 
on $\CX$, and choose a metric on $X_\sigma(\C)$ that induces the usual
topology. With $a\gg0$, for every $\theta \in X(\C_\sigma)$,

\item
\[ \underline{t}_1(\theta) \underline{t}_2(\theta) \geq 
   \underline{t}_1^a(\theta) \underline{t}_1^a(\theta) \geq1 \geq
   \underline{t}_1^a \bar{t}_2^a(\theta) \geq
   \underline{t}_1(\theta) \bar{t}_2(\theta). \]

\item
If $a$ is sufficiently big, for every $\theta \in X_\sigma$ 
\[ \limsup \left(-\inf_{F=(\alpha,f) \in \CF_D} \left\{
            \frac{D(F,\theta)}{a (\deg \alpha)^{1+\frac1t}} 
             \right\}     \right) \geq b > 0, \]
where $b$ only depends on $t$ and the degree and height of $\CX$; hence, 
the inequalities
\[ \underline{t}_1(\theta) \geq \underline{t}_1^a(\theta) \geq 
   \frac1{t(\theta)} \]
hold, and consequently 
\[ t(\theta) \geq \bar{t}_2^a(\theta) \geq \bar{t}_2(\theta). \]

\item
If $Y$ is any subscheme of $X$, then for all $\theta \in Y_\sigma(\C)$
with the exception of a subset of measure zero, the equalities
\[ t(\theta) =  \frac1{\underline{t}_1(\theta)} = 
   \frac1{\underline{t}_2^a(\theta)} =  
   \underline{t}_2(\theta) = \bar{t}_2(\theta) = 
   \underline{t}_2^a(\theta) = \bar{t}_2^a(\theta) \]
hold.

\end{enumerate}
\end{Satz}

The lower bounds for $\underline{t}_1$ and $\underline{t}_1^a$ in part two
are statements about the approximability of $\theta$ by algebraic points. 
It entails

\satz{Theorem} \label{con1}
In the situation of the previous Theorem, 
there exists a positive real number $b$ only depending
on $t(\theta)$ and the degree and height of the closure of $\psi(\theta)$, 
such that for any sufficiently big real number $a$, there is an infinite subset 
$M \subset \N$ such that for all $D \in M$ there exists an
algebraic point $\alpha_D \in \CX(\bar{\Q})$ fulfilling
\[ \deg(\alpha_D) \leq D^t, \quad h(\alpha_D) \leq a D^t, \quad
   \mbox{and} \quad \log |\alpha_D, \theta| \leq - ab D^{t+1}, \]
where $h(\alpha_D)$ denotes the height of $\alpha_D$, and
$|\cdot, \theta|$ the distance to the point $\theta$ with
respect to any metric on $X(\C_\sigma)$. 
\end{Satz}

The lower bound for $t(\theta)$ from Theorem \ref{approx}.2 may serve to
prove lower bounds for transcendence degrees and will be exploited in
the fifth paper of this series.

Theorem \ref{approx} will be proved in part 3 of this series (\cite{App3}).
The first inequality in its second part  
is best possible in the following sense: There is a constant
$B>0$ such that the set of points $\theta \in X(\C_\sigma)$ fulfilling the 
Theorem with $b$ replaced by $B$ has measure $0$ %(see \cite{mahler}).
For $t=1$, the Theorem in slightly different formulation was already 
proved in \cite{RW}.

The quantities $|X,\theta|$ and $|f_\theta|$ defining
$D(F,\theta)$ are important for applications, but rather ill-suited for proofs
as they don't fit into the framework of algebraic geometry: $|f_\theta|$
describes ``distances'' of points $\theta \in X_\sigma$ only to cycles of
codimension one, namely to $\di f$ for $f \in \Gamma(\CX,\CL)$, and
$|X,\theta|$ has no functorial properties and does not behave well
with respect to intersections.

To overcome these deficiencies, various authors (Philippon \cite{Ph},
Nesterenko \cite{Nes}) developed a new kind of distance between an 
effective cycles and a point in projective space that 
semi officially has been called algebraic distance.
Although this lead to proofs of hitherto unknown approximation 
Theorems (\cite{Ph1}), strong Theorems for approximation of transcendental 
points, that have good geometrical interpretation
and thus tie in well with the calculus of distances 
$|x,\theta|$, $|X,\theta|$ of points $\theta$ to points or cycles could
only be made in the one dimensional case (\cite{RW}).

This paper taking the concepts and results of \cite{RW} as a starting
point presents a geometrical approach to the concept of algebraic
distances in the framework of Arakelov geometry.
It allows to proof
amply stronger versions of arithmetic B\'ezout Theorems in the higher 
dimensional case, that accordingly have many new applications.
In the present context, the arithmetic B\'ezout Theorem will be used
to proof upper estimates for the algebraic distance of the intersection of 
two cycles
to a point $\theta$, given that both cycles have small algebraic distances
to $\theta$ compared with their degrees and heights, as well as
upper estimates for the distances of points $\alpha$ to $\theta$ when 
$\alpha$ belongs 
to a certain frame $(\alpha,f)$. Also, it will give upper bounds for
the approximational degrees as stated in the last two inequalities
of Theorem \ref{approx}.1.

Proving the existence of the points $\alpha$ claimed in the Theorems
above due to the metric B\'ezout Theorem reduces to
finding sufficiently many hyperplanes with small algebraic distance to a
point $\theta$ that intersect properly, hence produce cycles of higher
codimension, in particular points, with small algebraic distance to 
$\theta$. These hyperplanes can be found once one has good explicit upper 
bounds for arithmetic Hilbert
functions which will be proved in the second paper of this series.

\section{The main results}

Let $E = \Z^{t+1}$, and equip $E_\C:= E \otimes_{\Z} \C$
with the standard hermitian product. It
induces a hermitian metric on the line bundle $O(1)$ on the projective space
\[ \Pe^t = Proj(Sym(\check{E})), \]
which in turn defines an $L^2$-norm on the space of 
global sections $\Gamma(\Pe^t(\C),O(D))$, and a height $h(\CX)$ for any
effective cycle $\CX$ on $\Pe^t$.
A subscheme of $\Pe^t$ of dimension greater zero that has 
nonempty generic fibre will be called a subvariety.

Let further $|\cdot,\cdot|$ denote the Fubini-Study metric on $\Pe^t(\C)$, and
$\mu$ its K\"ahler form which also is the Chern form of $\overline{O(1)}$.
Finally for any projective subspace $\Pe(W) \subset \Pe^t(\C)$, denote
by $\rho_{\Pe(W)}$ the function
\[ \rho_{\Pe(W)}:
   \Pe^t(\C) \setminus \Pe(W) \to \R, \quad x \mapsto \log |x,\Pe(W)|. \]

For $\CX \in Z^p(\Pe^t)$ an effective cycle of pure codimension $p$, and
$\theta$ a point in  $\Pe^t(\C) \setminus supp(X_\C)$  
the logarithm of the distance $\log |\theta,X|$ is defined to be the
minimum of the restriction of $\rho_{\theta}$ to $X$. 
There are various different 
definitions of the algebraic distance of $\theta$ to $X = \CX(\C)$ 
all identical modulo a constant times $\deg X$; the simplest being

\satz{Definition}
With the above notations, and $\Lambda_{\Pe(W)}$ the Levine 
form of a projective subspace $\Pe(W)$, define the algebraic distance of 
$\theta \not\in supp(X_\C)$ to $X$ as
\[ D(\theta,X) := {\mbox{sup} \atop \Pe(W)}
   \int_{X(\C)} \Lambda_{\Pe(W)} -
   \deg X \sum_{n=1}^q \sum_{m=0}^{t-q} \frac{1}{m+n}, \]
where the supremum is taken over all spaces $\Pe(W) \subset \Pe^t$ of 
codimension $t+1-p$ that contain $\theta$.
In case $supp(X)=\emptyset$,
the algebraic distance $D(\theta,X)$ is defined to be zero. 

%This distance is always nonpositive. 
%Finally the canonical linear form $\varphi_\theta (X)$ is defined as
%\[ \varphi_\theta(X) := \Sigma_\theta(X) + h(X). \]
\end{Satz}

\satz{Proposition} \label{haupt1}
Let $X \in Z_{eff}(\Pe^t_\C)$ be an effective cycle and $\theta \in \Pe^t(\C)$
a point that is not contained in the support of $X$. 
There are effectively computable constants $c,c'$ only depending
on $t$ such that
\[ \deg(X) \log |\theta,X(\C)| \leq D(\theta,X) + c \deg X \leq 
   \log |\theta,X(\C)| + c' \deg X. \]
\end{Satz}

\satz{Proposition} \label{haupt2}
Let $f \in \Gamma(\Pe^t_\Z,O(D))$, and $\CX = \di f$. Then,
\[ h(\CX) \leq \log |f_D|_{L^2} + D \sigma_t, \quad \mbox{and} \]
\[ D(\theta,X) + h(\CX) = \log |f_\theta| + D \sigma_t, \]
where the $\sigma_i's$ are certain constants.
\end{Satz}

\satz{Theorem} {\bf (Metric B\'ezout Theorem}) \label{bezout} 
Let $p+q \leq t+1$, and
$\CX,\CY \in \Pe^t_\Z$ be effective cycles of pure codimensions $p$ and $q$ 
respectively intersecting properly. 
There is an effectively computable positive constant $e$,
only depending on $t,p,q$, and for each
$\theta$ a point in $\Pe^t(\C) \setminus (supp(X_\C\cup Y_\C))$ there exists
a map 
\[ f_{X,Y}: \underline{\deg X + \deg Y} 
  \to \underline{\deg X} \times \underline{\deg Y} \]
from the set of natural numbers less or equal
$\deg X + \deg Y$ to the set of natural numbers
less or equal $\deg X$ times the set of natural numbers less or equal
$\deg Y$ such that the maps $pr_1 \circ f_{X,Y}, pr_2\circ f_{X,Y}$
are monotonously increasing and surjective, $f_{X,Y}$ is a right inverse to
the sum and for every
$T\in \underline{\deg X+\deg Y}$, with $(\nu,\kappa) = f_{X,Y}(T)$, 
the inequality
\[ \nu \kappa \log |\theta,X+Y| + D(\theta,X.Y) + h(\CX . \CY) \leq \]
\[  \kappa D(\theta,X) + \nu D(\theta,Y) +
   \deg Y h(\CX) + \deg X h(\CY) + e \deg X \deg Y \]
holds.
\end{Satz}

For $t=1$, this Theorem has been proved in \cite{RW}. 

\satz{Corollary}  \label{haupt3}
In the situation of the Theorem,
\begin{enumerate}
\item
if either $D(X,\theta) \leq \log |Y,\theta|$ or $|X,\theta| \leq |Y,\theta|$, 
then with $e'$ an effectively computable constant only depending on $t,p,q$, 
\[ D(\theta,X.Y) + h(\CX . \CY) \leq D(\theta,Y) +
   \deg Y h(\CX) + \deg X h(\CY) + e' \deg X \deg Y. \]
\item
in any case
\[ D(\theta,X.Y) + h(\CX . \CY) \leq \max(D(\theta,X),D(\theta,Y)) +
   \deg Y h(\CX) + \deg X h(\CY) + \]
\[ e' \deg X \deg Y. \]
\end{enumerate}
\end{Satz}

A variant of this corollary
has been proved in \cite{Ph} in case $\CX$ is a hypersurface.

The crucial idea for proving the Theorem is to express
the algebraic distance of an effective cycle $X$ of pure codimension $p$
to a point $\theta$ as the sum of the
distances of $\theta$ to certain points lying on $X$, namely the points
forming the intersection of $X$ with a suitable projective subspace of
$\Pe^t$ of dimension $p = \mbox{codim} \; \CX$. More precisely,

\satz{Theorem}
For $p \leq t$ there are effectively computable
constants $c,c'$ only depending on $p$ and $t$,
such that for all effective cycles $X$ of pure codimension
$p$ in $\Pe^t_\C$, and points $\theta \in \Pe^t(\C)$ 
not contained in $supp(X)$, there is a subspace $\Pe(F) \subset \Pe^t_\C$ of 
codimension $t-p$ containing $\theta$, and properly intersecting $X$, such that
\[ D(\theta,X) \leq 
   \sum_{x \in supp(X . \Pe(F))} n_x \log |x,\theta| + c \deg X \leq
   D(\theta,X) + c' \deg X, \]
where the $n_x$ are the intersection multiplicities of $X$ and $\Pe(F)$ at $x$.
\end{Satz}
This Theorem is proved in section 4.

\section{Arakelov varieties}

This section mainly collects various well known facts about Arakelov varieties,
most of which can be found in \cite{SABK}, \cite{GS1}, or \cite{BGS}.

Let $\CX$ a regular, flat, projective scheme of relative dimension $d$ over 
$\spec \; \Z$. Such a scheme is called a projective arithmetic variety over
$\spec \; \Z$
The base extensions of $\CX$ to $\Q$ and $\C$, as well as their $\C$ valued
point $X(\C)$ will all be denoted by $X$ if no confusion arises.
On the $\C$- valued points $X(\C)$ we have the space of
smooth forms of type $(p,p)$ invariant under complex conjugation
$F_\infty$ denoted by $A^{p,p}$, the space 
$\tilde{A}^{p,p} := A^{p,p}/(\mbox{Im} \partial + \mbox{Im} \bar{\partial})$,
and the space of currents $D^{p,p}$ which is the space of Schwartz continuous
linear functionals on $A^{d-p,d-p}$, and 
$\tilde{D}^{p,p} = D^{p,p}/(\mbox{Im} \partial + \mbox{Im} \bar{\partial})$. 
On $D^{p,p}$ the maps 
$\partial, \bar{\partial}, d = \partial + \bar{\partial}, d^c = \partial-
\bar{\partial}$ are defined by duality. 

A cycle $Y \subset Z^p (X_\C)$ of pure codimension $p$ defines a 
current $\delta_Y \in D^{p,p}$ by
\[ \omega \in A^{d-p,d-p} \mapsto \sum_i n_i \int_{Y_i} \omega, \]
where $Y = \sum_i n_i Y_i$ is the decompositions into irreducible 
components leading an embedding
$\iota: A^{p,p} \hookrightarrow D^{p,p}, \omega \mapsto [\omega]$.
The integrals are defined by resolution of singularities, see \cite{GS1}
or \cite{SABK}.
A Green current $g_Y$ for $Y$ is a current of type $(p-1,p-1)$ such that
\[ d d^c g_Y + \delta_Y \in \iota(A^{p,p}). \]
A densely defined form $g_Y$ on $X$ is called a Green form for $Y$ if
$[g_Y]:=\iota(g_Y)$ exists and is a Green current for $Y$; 
it is called of logarithmic type along
$Y$ if it has only logarithmic singularities at $Y$ 
(see \cite{SABK} Def.\@ II.\@2.\@3).

If $y$ is a point in $\CX^{(p-1)}$, that is $\CY = \overline{\{y\}}$ is a 
closed integral sub scheme of codimension $p$ in $\CX$, a rational function
$f \in k(y)^*$ gives rise to the Green form of log type
$-\log |f|^2$ for $\di(f)$. (\cite{SABK}, III.\@1)

The group of arithmetic cycles $\hat{Z}^p(\CX)$ consisting of the
pairs $(\CY,g_Y)$ where $\CY$ is a cycle of pure codimension $p$ and 
$g_Y$ a Green current for $Y(\C)$, thus contains the subgroup
$\hat{R}^p(\CZ)$ generated by the pairs $(\di(f),- [\log|f|^2])$, with
$f \in k(y)^*, y \in \CX^{(p-1)}$ and pairs 
$(0,\partial(u) +\bar{\partial}(v))$,
where $u$ and $v$ are currents of type $(p-2,p-1)$ and $(p-1,p-2)$ 
respectively. If we set 
$\tilde{Z}^p(\CX) :=  \hat{Z}^p(\CX)/
(\mbox{Im} \partial + \mbox{Im}\bar{\partial})$,
and
$\tilde{R}^p(\CX) :=  \hat{R}^p(\CX)/
(\mbox{Im} \partial + \mbox{Im}\bar{\partial})$, we get

\satz{Definition}
The arithmetic Chow group $\widehat{CH}^p(\CX)$ of codimension $p$  of $\CX$
is defined as the quotient $\tilde{Z}^p(\CX)/ \tilde{R}^p(\CX)$.
\end{Satz}

\satz{Examples}
\begin{enumerate}
\item
For $S = \spec \Z$, it is easily calculated (see \cite{BGS}, 2.\@ 1.\ 3)
\[ \widehat{CH}^0(S) \cong \Z,  \quad \widehat{CH}^1(S) \cong \R, \quad
   \widehat{CH}^p(S) = 0, \quad \mbox{if} \quad p>1. \]
The isomorphism of $\widehat{CH}^1(S)$ to $\R$ is usually denoted
$\widehat{deg}$.

\item
If $\bar{\CL}$ is an ample metric line bundle on $\CX$, and
$f \in \Gamma(\CX,\CL)$ is a global section then
$[-\log |f|^2]$, by the Poincar\'e-Lelong formula, is a Green current
for $\di \; f$; we have 
$d d^c [-\log |f|^2] + \delta_{\di f} = c_1(\bar{L})$ the
Chern form of $\bar{L}$. The class 
\[ \hat{c}_1 (\bar{\CL}) := (\di f,[-\log|f|^2]) \in \widehat{CH}^1(\CX) \]
is called the first Chern class of $\bar{\CL}$.
\end{enumerate}
\end{Satz}

For the purposes of this paper a subgroup of the arithmetic Chow
group, the Arakelov Chow group will play a much more important role.
Suppose $X(\C)$ is equipped with a K\"ahler metric with K\"ahler form
$\mu$. The pair $(\CX,\mu)$ is denoted $\bar{\CX}$, and called
an Arakelov variety. If
$H^{p,p}(X)$ denotes the harmonic forms with respect 
to the chosen metric, and 
\[ H: D^{p,p}(X) \to H^{p,p}(X) \]
the harmonic projection, define 
\[ Z^p(\bar{\CX}) := \{ (\CY,g_Y) | d d^c g_Y + \delta_Y \in H^{p,p}(X) \}
   \subset \hat{Z}^p(X). \]
As $\hat{R}^p(\CX) \subset Z^p(\bar{\CX})$, one can define the
Arakelov Chow group  
\[ CH^p(\bar{\CX}) := Z^p(\bar{\CX})/\hat{R}^p(\CX). \]

For $[(\CY,g_Y)] \in CH^p(\bar{\CX})$ we have
\begin{equation} \label{hpz}
d d^c g_Y + \delta_Y = \omega_Y = H(\omega_Y) = 
H(d d^c g_Y) + H(\delta_Y)   = H(\delta_Y).
\end{equation}
A Green $g_Y$ current with $d d^c g_Y + \delta_Y = H(\delta_Y)$ 
i.\@ e.\@ $(\CX,g_X) \in Z(\bar{\CX})$ is
called admissible.

In \cite{GS1} Gillet and Soul\'e define the star product
\[ [g_Y] * g_Z := [g_Y] \wedge \delta_Z + g_Z \wedge \omega_Y \]
of a Green current $[g_Y]$ coming from a Green form $g_Y$ of log type
along $Y$ with a Green current $g_Z$ for properly intersecting $Y$, and $Z$.

\satz{Proposition} \label{schprodukt}
If $d d^c [g_Y] + \delta_Y = [\omega_Y]$, and 
$d d^c g_Z + \delta_Z = [\omega_Z]$, then
\[ d d^c([g_Y] * g_Z) + \delta_{Y . Z} = [\omega_Y \wedge \omega_Z]. \]
Further, the star product is commutative and associative modulo 
$\mbox{Im} \partial + \mbox{Im} \bar{\partial}$. As by \cite{GS1} every
cycle has a Green form of log type, there is an intersection product
\[ \widehat{CH}^p(\CX)\times \widehat{CH}^q(\CX) \to \widehat{CH}^{p+q}(\CX), 
   \quad ([\CY,[g_Y]],[\CZ,g_Z]) \mapsto ([\CY . \CZ], [g_Y] * g_Z), \]
where $\CY,\CZ$ are chosen to intersect properly, which is commutative and 
associative.

If on $\bar{\CX}$ the product of two harmonic forms is always harmonic,
the above product makes $CH(\bar{\CX})^*$ into a subring of 
$\widehat{CH}^*(\CX)$.
\end{Satz}

\proof \cite{GS1}, II, Theorem 4.

\satz{Proposition} \label{funk}
Let $\CX,\CY$ be arithmetic varieties over $\spec \Z$, and
$f: \CX \to \CY$ a morphism.

\item
For $(Z,g_Z) \in \hat{Z}^p(Y)$ we have 
$d d^c f^*g_Z + \delta_{f^*(Z)} = f^* \omega_Z$, and the densely defined map
\[ f^* : \hat{Z}^p(Y) \to \hat{Z}^p(X), \quad
   (Z,g_Z) \mapsto (f^*(Z), f^*g_Z), \]
induces a multiplicative pull-back homomorphism 
$f^*: \widehat{CH}^p (\CY) \to \widehat{CH}^p(\CX)$.

\item
If $f$ is proper, $f_\Q: X_\Q \to Y_\Q$ is smooth, and $X,Y$ are
equidimensional, then  
$d d^c f_* g_Z + \delta_{f_* Z} = f_* \omega_Z$ for any
$(\CZ,g_Z) \in \hat{Z}^p(\CX)$. This induces a push-forward homomorphism
\[ f_*: \widehat{CH}^p(\CX) \to \widehat{CH}^{p-\delta}(\CY), \quad
   (\delta:= \mbox{dim} Y - \mbox{dim} Z). \]

If $f_* (f^* respectively)$ map harmonic forms to harmonic forms,
they induce homomorhpisms of the Arakelov Chow groups.

\end{Satz}

\proof
\cite{SABK}, Theorem III.\@ 3.

The following Proposition enables calculations in
$\widehat{CH}^*(\CX)$ and $CH^*(\bar{\CX})$.

\satz{Proposition} \label{sequ}
Let $a: A^{p-1,p-1}(X) \to \widehat{Ch}^p(\CX)$ be the map
$\eta \mapsto [(0,\eta)]$, and $\zeta: \widehat{Ch}^p(\CX) \to Ch(\CX)$ the
map $[(\CY,g_Y)] \mapsto [\CY]$.  
With $\tilde{Z}^p(\bar{\CX}) = Z^p(\bar{\CX})/
(\mbox{Im} \partial + \mbox{Im} \bar{\partial})$,
the diagram

\vspace{6mm}

\hspace{-2cm}
\unitlength=1.00mm
\special{em:linewidth 0.4pt}
\linethickness{0.4pt}
\begin{picture}(160.00,67.00)
%\put(17.00,10.00){\vector(1,0){12.00}}
\put(61.00,13.00){\makebox(0,0)[cc]{$a$}}
\put(61.00,38.00){\makebox(0,0)[cc]{$a$}}
\put(61.00,63.00){\makebox(0,0)[cc]{$a$}}
\put(95.00,13.00){\makebox(0,0)[cc]{$\zeta$}}
\put(95.00,38.00){\makebox(0,0)[cc]{$\zeta$}}
\put(95.00,63.00){\makebox(0,0)[cc]{$\zeta$}}
\put(81.00,48.00){\makebox(0,0)[cc]{$pr$}}
\put(116.00,48.00){\makebox(0,0)[cc]{$pr$}}
\put(80.00,23.00){\makebox(0,0)[cc]{$\iota$}}
\put(45.00,23.00){\makebox(0,0)[cc]{$\iota$}}
\put(55.00,10.00){\vector(1,0){12.00}}
\put(42.00,30.00){\vector(0,-1){15.00}}
\put(77.00,30.00){\vector(0,-1){15.00}}
%\put(17.00,45.00){\vector(1,0){12.00}}
\put(55.00,35.00){\vector(1,0){12.00}}
\put(87.00,35.00){\vector(1,0){14.00}}
\put(122.00,35.00){\vector(1,0){20.00}}
\put(42.00,55.00){\vector(0,-1){15.00}}
\put(77.00,55.00){\vector(0,-1){15.00}}
\put(112.00,55.00){\vector(0,-1){15.00}}
\put(112.00,30.00){\vector(0,-1){15.00}}
%\put(17.00,80.00){\vector(1,0){12.00}}
\put(55.00,60.00){\vector(1,0){12.00}}
\put(85.00,60.00){\vector(1,0){18.00}}
\put(122.00,60.00){\vector(1,0){20.00}}
%\put(12.00,10.00){\makebox(0,0)[cc]{$0$}}
%\put(12.00,45.00){\makebox(0,0)[cc]{$0$}}
%\put(12.00,80.00){\makebox(0,0)[cc]{$0$}}
\put(147.00,60.00){\makebox(0,0)[cc]{$0$}}
\put(147.00,35.00){\makebox(0,0)[cc]{$0$}}
\put(147.00,10.00){\makebox(0,0)[cc]{$0$}}
\put(42.00,10.00){\makebox(0,0)[cc]{$\tilde{A}^{p-1,p-1}(X)$}}
\put(42.00,35.00){\makebox(0,0)[cc]{$H^{p-1,p-1}(X)$}}
\put(42.00,60.00){\makebox(0,0)[cc]{$H^{p-1,p-1}(X)$}}
\put(122.00,10.00){\vector(1,0){20.00}}
\put(77.00,10.00){\makebox(0,0)[cc]{$\widehat{CH}^p(\CX)$}}
\put(77.00,35.00){\makebox(0,0)[cc]{$CH^p(\bar{\CX})$}}
\put(77.00,60.00){\makebox(0,0)[cc]{$\tilde{Z}^p(\bar{\CX})$}}
\put(112.00,10.00){\makebox(0,0)[cc]{$CH^p(\CX)$}}
\put(87.00,10.00){\vector(1,0){14.00}}
\put(112.00,35.00){\makebox(0,0)[cc]{$CH^p(\CX)$}}
\put(112.00,60.00){\makebox(0,0)[cc]{$Z^p(\CX)$}}
\end{picture}

is commutative, and the rows are exact.
\end{Satz}

\proof
\cite{GS1}

\vspace{3mm}

If $\CY \in Z^p(\CX)$ and $g_Y,g'_Y$ are two admissible Green currents for
$Y$, the exactness of the first row implies 
$g_Y-g'_Y = \eta \in H^{p-1,p-1}(X)$. Hence, the projection
of $g_Y$ to the orthogonal Complement of $H^{p-1,p-1}$ in
$\tilde{D}^{p-1,p-1}$ is independent of $g_Y$, and one can define
the map $s:Z^p(X) \to Z^p(\bar{X})$ by
$Y \mapsto (Y,g_Y)$ where $g_Y$ is the unique Green form of log type for $Y$ 
which is orthogonal to $H^{p-1,p-1}(X)$. 
Then, $s$ defines a splitting of the first
exact sequence, and induces a pairing
\begin{equation} \label{hpaarung}
\widehat{CH}^p(\CX) \times Z^q(\CX) \to \widehat{CH}^{p+q}(\CX), \quad
(y,\CZ) \mapsto (y|\CZ) = y . [s(\CZ)].
\end{equation}
A Green current $g_Y$ with $(\CY,g_Y)=s(\CY)$ i.\@ e.\@ 
$(\CY,g_Y) \in Z^p(\bar{\CX})$, and $g_Y$ is orthogonal to $H^{d-p,d-p}$
is called ($\mu$-)normalized. It is unique modulo
$\mbox{Im} \partial + \mbox{Im} \bar{\partial}$.

\satz{Definition}  \label{hoehe}
For $\pi:\CX \to \spec \; \Z$ a projective arithmetic variety,
$\bar{\CL}$ a metrized ample line bundle on $\CX$ with Chern form
$c_1(\bar{L})$ equal to $\mu$, 
the height of an effective cycle $\CZ \in Z^p(\CX)$ is defined as
\[ h(\CZ) := \widehat{\deg} (\pi_*(\hat{c}_1(\bar{\CL})^{d-p+1}|\CZ)) \in \R. \]
If $\hat{c}_1^{d+1-p}(\bar{\CL}) = [(\CY,g_Y)]$, and
$supp(\CY) \cap supp(\CZ) = 0$, this is equal to
\[ \frac12 \int_{Z(\C)} g_Y. \]
\end{Satz}

\satz{Proposition} \label{hfunk}
Let $\CX,\CY$ be regular, projective, flat schemes over $\spec \Z$, and
$f: \CX \to \CY$ a morphism. Further, let 
$p,q$ be natural numbers with $p+q = d+1=\dim \CX$, and
$(\CZ,g_Z) \in \widehat{CH}^p(\CX), \\ \CW \in Z^q(\CY)$. 
If $\dim f(\CW) = \dim(\CW)$, we have
\[ \left( f^*(\CZ,g_Z)| \CW \right) = \left( (\CZ,g_Z)| f_*(\CW) \right). \]

If $f$ is flat, and surjective, has smooth restriction to every 
component of $X_\Q$, 
and $\CX, \CY$ have constant dimension, then with
$\delta = \dim \CX - \dim \CY$,
$(\CZ,g_Z) \in \widehat{CH}^p(\CX)$, $\CW \in Z^{d+1-p-\delta} (\CY)$, we have
\[ \left( (\CZ,g_Z) | f^*(\CW) \right) = \left( f_*([\CZ,g_Z])|\CW \right). \]
\end{Satz}

\proof
\cite{BGS}, Proposition 2.\@3.\@1, (iv),(v).

\satz{Proposition} \label{hsch}
\begin{enumerate}

\item
Let $\CY$ be an effective cycle of pure codimension $p$ on $\CX$,
$\CL$ an ample line bundle on $\CX$, and
$f$ a global section of $\CL^{\otimes D}$ on $\CX$, whose
restriction to $Y$ is nonzero. Then,
\[ h(\CY . \di f) = 
   D h(\CY) + \int_{X(\C)} \log ||f|| \mu^{d-p} \delta_Y. \]

\item
Assume that on the variety $\bar{X}$ the product of two harmonic forms
is always harmonic, and $\CY,\CZ$ are
effective cycles of pure codimensions $p$ and $q$ respectively,
intersecting properly. With $s(\CY) = (\CY,g_Y)$,
\[ s(\CY . \CZ) = s(\CY) . s(\CZ)  
   -\frac12 a(H(g_Y \delta_Z c_1(\mu^{d+1.p.q}))), \]
and consequently,
\[ h(\CY . \CZ) =
   \pi_*(\hat{c}_1(\CL)^{d+1-p-q} 
   . s(\CY) . s(\CZ) + a(H(g_Y \delta_Z))). \]
\end{enumerate}
\end{Satz}

\proof
1. \cite{BGS}, Proposition 3.2.1 (iv).

%1. Let $\hat{c}_1(\bar{\CL})^{d-p} = (\CZ, g)$, with
%$[\CZ] = c_1{\CL}^{d-p}$, and $s(\CY) = (\CY,g_Y)$. 
%Then,
%\begin{eqnarray*}
%D (\hat{c}_1^{d+1-p}|\CY) &=& D \hat{c}_1^{d+1-p} . (\CY,g_Y) =
%   \hat{c_1}^{d-p} . (\di f, -\log |f|^2) . (\CY,g_Y)  \\
% &=&
%(c_1(\CL)^{d-p} . \di f, g \delta_{\di f} - \log |f|^2 \mu^{d-p}) 
%   . (\CY,g_Y) \\
%&=& (c_1(\CL)^{d-p} . \di f . \CY, 
%    g \delta_{\di f} \delta_Y - \log |f|^2 \mu^{d-p} \delta_Y + 
%    g_Y \mu^{d+1-p}). 
%\end{eqnarray*}
%On the other hand
%\begin{eqnarray*}
%(\hat{c}_1(\CL)^{d-p}|\CY . \di f) &=& 
%   \hat{c}_1^{d-p} . s(\CY . \di f)  \\
%       &=&
%    (c_1(\CL)^{d-p},g) . 
%   (\di f . Y, -\log |f|^2 \delta_Y + D g_Y \mu \\ 
%    & & -H(-\log |f|^2 \delta_Y + D g_Y \mu)) \\
%      &=&  (c_1(\CL)^{d-p} . \di f . Y, 
%    g \delta_{\di f} \delta_Y - \log |f|^2 \mu^{d-p} \delta_Y + 
%    g_Y \mu^{d+1-p}) \\
%    & & -a(H(-\log |f|^2 \delta_Y \mu^{d-p}+ D g_Y \mu^{d+1-p})). 
%\end{eqnarray*}
%Hence,
%\[ D (\hat{c}_1^{d+1-p}|Y)- (\hat{c}_1(\CL)^{d-p}|Y . \di f) =
%   -a(H((-\log |f|^2 \delta_Y + D g_Y \mu) \mu^{d-p})). \] 
%Consequently,
%\begin{eqnarray*}
%D h(Y) - h(Y . \di f) &=&
%   \pi_*(H(-\log |f|^2 \delta_Y \mu^{d-p}+ D g_Y \mu^{d+1-p})) \\
%&=& \int_{Y(\C)} -\log |f|^2 \mu^{d-p} + D \int_{X(\C)} g_Y \mu^{d+1-p} \\
%&=& \int_{Y(\C)} -\log |f|^2 \mu^{d-p}, 
%\end{eqnarray*}
%since $g_Y$ is orthogonal to the harmonic form $\mu^{d+1-p}$.

\vspace{2mm}

2.
Let $s(\CY) = (\CY,g_Y)$, and $s(\CZ) = (\CZ,g_Z)$. Then,
\[ s(\CY) . s(\CZ) = (\CY . \CZ, g_Y \delta_Z + H(\delta_Y) g_Z). \]
As the form
\[ d d^c(g_Y \delta_Z + H(\delta_Y) g_Z) + \delta_{Y . Z} =
   H(\delta_Y) H(\delta_Z) \]
is harmonic by assumption, 
$(\CY . Z, g_Y \delta_Z + H(\delta_Y) g_Z) \in Z^{p+q} (\bar{X})$, and
\[ s(\CY . \CZ) = 
   (\CY . \CZ, g_Y \delta_Z + H(\delta_Y) g_Z- 
               H(g_Y \delta_Z + H(\delta_Y) g_Z)). \]
Since multiplication with a harmonic forms leaves the
space of harmonic forms invariant, it also leaves the space 
of forms orthogonal to the harmonic forms invariant; hence
$H(\delta_Y) g_Z$ is orthogonal to the space of harmonic forms,
and $H(H(\delta_Y) g_Z) = 0$, implying
\[ s(\CY) . s(\CZ) = s(\CY . \CZ) + a(H(g_Y \delta_Z)). \]
The claim about the heights follows by multiplying the last equality
with \\ $\hat{c}_1(\bar{\CL})^{d-p-q+1}$ and applying $\pi_*$.

%\satz{Remark}
%{\bf 2 besser als 1, weil...}
%\end{Satz}

An important tool for making estimates is the concept of positive Green
forms: A smooth form $\eta$ of type $(p,p)$ on a complex manifold 
is called positive if for any complex sub manifold 
$\iota:V \to X$ of dimension $p$, the volume form
$\iota^* g_Y$ on $V$ is nonnegative, i.\@ e.\@ for each point
$v \in V$, the local form $(\varphi^* g_Y)_v$ is either zero or induces the 
canonical local orientation at $v$.

\satz{Lemma} \label{prodpos}
Let $X,Y$ be complex manifolds, and $\eta$ a positive form
of type $(p,p)$ on $X$.
\begin{enumerate}

\item
For any holomorphic map $f: Y \to X$, the form $f^* \eta$ is positive. 

\item
If $g: X\to Y$ is a smooth holomorphic map whose restriction to the support of
$\eta$ is proper, the form $g_*\eta$ is positive.

\item
For any positive form $\omega$ of type $(1,1)$ 
the form $\omega \wedge \eta$ is positive.

\end{enumerate}
\end{Satz}

\proof
\cite{BGS}, Proposition 1.\@1.\@4.

\subsection{Projective Space}

Let $M$ be a free $\Z$ module of rank $t+1$, and
$\Pe^t = Proj(Sym(\check{M}))$ the projective space with
structural morphism
$\pi: \Pe^t \to \spec \Z$. If 
$M_\C = \C \otimes \Z^{t+1}$ is equipped with a hermitian product, this
induces a metric on the line bundle
$O(1)$ on $\Pe^t_\C$ and the Fubini-Study metric on $\Pe^t(\C)$.
The Chern form $\mu = c_1(\overline{O(1)})$ equals the K\"ahler form 
corresponding to this metric. 

For any torsion free submodule $N \subset M$ define 
$\widehat{\deg}(N \otimes_\Z \Q) = \widehat{\deg}(N)$ as minus the logarithm
of the covolume of $N$ in $N_\C = N \otimes_\Z \R$.
We will always assume that the hermitian product chosen in such a way that
$\widehat{\deg}(M) = 0$.

%This induces   and the vector space 
%$E_\C = E \otimes_\Z \C$ be equipped with a hermitian product. The module
%$E$ together with this product forms an arithmetic module $\bar{E}$. For 
%any $F \subset E$ define the arithmetic 
%$\widehat{\deg} (\bar{F})$ as minus the 
%logarithm of the covolume of the lattice $F$ inside $F_\C$ with respect
%to the hermitian product that comes from restricting the hermitian product
%on $E_\C$ to $F_\C$.  
%For simplicitiy we will always take $E=\Z^{t+1}$ equipped with the
%standard hermitian product. In this case the degree of any submodule $F$
%will always be at most $0$.
%Denote the dual space of $E$ by $\check{E}$. It canonically carries
%a hermitian product, and the corresponding arithmetic module will be denoted
%by $\check{\bar{E}}$.

%Let $\Pe^t = \Pe(E) = Proj(\sym(\check{E}))$, and 
%$\pi: \Pe^t \to \spec(\Z)$ the structural morphism. On $\Pe^t$ there
%is a short exact sequence of vector bundles
%\[ 0 \to Q  \to \pi^*(\check{E}) \to O(1) \to 0, \]
%where $Q_p$ is the kernel of the map that evaluates $f \in \check{E}$ at
%the point $p \in \Pe(E)$. As $\pi^*(\check{E})$ canonically carries
%a smooth hermitian metric, the line bundle $O(1)$ may be endowed with
%the quotient metric thus giving rise to a metrized line bundle 
%$\overline{O(1)}$.

Then, the height of a projective subspace $\Pe(F) \subset \Pe^t$ of
dimension $p$ equals
\begin{equation} \label{hpr}
h(\Pe(E)) = - \widehat{\deg}(\bar{F}) + \sigma_p,
\end{equation}
where the number
\[ \sigma_p := \frac{1}{2} \sum_{k=1}^p \sum_{m=1}^k \frac{1}{m} \]
is called the $p$th Stoll number (\cite{BGS}, Lemma 3.\@3.\@1).
The height of any effective cycle $\CZ \in Z^*(\Pe^t)$ is nonnegative
(\cite{BGS}, Proposition 3.\@2.\@4).

The space of harmonic forms $H^{p,p}(\Pe^t)$ with respect to the chosen metric
is one dimensional with
generator $\mu^p$. By Proposition \ref{sequ}, together with the 
definition of the map $s$, 
an element in $CH^p(\bar{\Pe}^t)$ may be written as
$\alpha \hat{\mu}^p + \beta a(\mu^p)$ with 
$\hat{\mu} := \hat{c}_1(\overline{O(1)})$, $\alpha \in \Z, \beta \in \R$.
 As $\zeta \circ s = id$, the degree of any
$\alpha \hat{\mu}^p + \beta a(\mu^p)$ equals $\alpha$.
One easily calculates $a(\mu^p) . a(\mu^q) = 0$, and 
$\hat{\mu}^p a(\mu^q) = a(\mu^{p+q})$.
This, together with Proposition \ref{hsch}.\@2 
implies the Proposition (\cite{BGS}, 5.\@ 4.\@3)

\satz{Proposition} \label{vorbezout}
Let $\CX,\CY$ be effective cycles of pure codimension $p$ and $q$ respectively
in $\Pe^t$ intersecting properly. With $(\CX,g_X) = s(\CX)$,
\begin{enumerate}
\item
\[ h(\CX . \CY) = \deg X h(\CY) + \deg Y h(\CX) -
   \frac12 \int_{\Pe^t} g_X \delta_Z \mu^{t+1-p-q}\; - \;  
   \sigma_t \deg X \deg Y. \] 
\item
With $c_1(p,q,t) := \sigma_{p+q-t}+\sigma_t-\sigma_p-\sigma_q$ and 
$\bar{c}_1(p,q,t) = c_1(p,q,t) + \frac{p+q-t-1}2$, the inequality
\[ -\frac12 \int_{\Pe^t} g_X \delta_Z \mu^{t+1-p-q} \leq c_1 \deg X \deg Y \]
holds.

\end{enumerate}
\end{Satz}

\proof
1. Assume 
\[ [s(X)] = \alpha \hat{\mu^p} + \beta a(\mu^p), \quad \mbox{and} \quad
   [s(Y)] = \alpha' \hat{\mu^q} + \beta' a(\mu^q). \]
Then $\alpha = \deg X, \alpha' = \deg Y$, and
\[ h(X) = \pi_*(\hat{\mu}^{t+1-p} (\alpha \hat{\mu}^p + \beta \mu^p)) = 
   \alpha \pi_*(\hat{\mu}^{t+1}) + \beta \pi_*(a(\mu^{t+1})) =
   \alpha \sigma_t + \beta. \]
Similarly $h(Y) =  \alpha' \sigma_t + \beta'$. Next, by Proposition
\ref{hsch}.\@2,
\[ h(X . Y) = 
   \pi_*(\hat{\mu}^{t+1-p-q} (\alpha \hat{\mu}^p + \beta \mu^p)
                            (\alpha' \hat{\mu}^p + \beta' \mu^p)) 
   -\frac12 \int_{\Pe^t} g_Y \delta_Z \mu^{t+1-p-q}. \]
The proposition thus follows from the calculation
\[ \hat{\mu}^{t+1-p-q} (\alpha \hat{\mu}^p + \beta a(\mu^p))
                            (\alpha' \hat{\mu}^p + \beta' a(\mu^p)) = \]
\[ \hat{\mu}^{t+1-p-q} (\alpha \alpha' \hat{\mu}^{p+q} +
   \alpha \beta' \hat{\mu}^p a(\mu^q) + \alpha' \beta \hat{\mu}^q a(\mu^p)) =
   \alpha \alpha' \hat{\mu}^{t+1} +
   (\alpha \beta' + \alpha' \beta) a(\mu^{t+1}), \]
and thus 
\begin{eqnarray*}
 h(\CX . \CY) &=& 
   \alpha \alpha' \sigma_t+ \alpha \beta'+\alpha' \beta 
   -\frac12 \int_{\Pe^t} g_Y \delta_Z \mu^{t+1-p-q} \\
  &=& \deg X h(\CY) + \deg Y h(\CX) - \alpha \alpha' \sigma_t 
   -\frac12 \int_{\Pe^t} g_Y \delta_Z \mu^{t+1-p-q}. 
\end{eqnarray*}

\vspace{2mm}

2. \cite{BGS}, Proposition 5.1.1 together with the proof of Theorem
5.4.4,(ii).

\vspace{2mm}

Let $F_\C \subset M_\C$ be a sub vector space of codimension $p$ 
with orthogonal complement $F_\C^\bot$, and $pr_{F^\bot}$ the projection
to the orthogonal complement.
Then, on $M_\C \setminus \{0 \}$ (resp.\@ $M_\C \setminus F_\C$)
the functions $\rho(x)=\log |x|^2$ 
(resp.\@ $\tau(x) = \log |pr_{F^\bot}(x)|^2$) are
defined, and give rise to the $(1,1)$-forms 
$\mu_M := d d^c \rho$ on $\Pe(E_\C)$, and 
$\lambda_{M,F} := d d^c \tau$ on $\Pe(M_\C) \setminus \Pe(F_\C)$ and 
a function $\rho -\tau$ on $\Pe(M_\C) \setminus \Pe(F_\C)$. Here,
$\mu_M$ is just the Chern form of the metrized line bundle $\overline{O(1)}$,
and $(\rho-\tau)(x)$ is $-\frac12$ times the logarithm of the Fubini-Study
distance of $x$ to $\Pe(F_\C)$.
With these notations, the so called Levine form
\begin{equation} \label{Levine}
\Lambda_{\Pe(F)} := (\rho-\tau) \sum_{i+j = p-1} \mu_M^i \lambda_{M,F}^j
\end{equation}
is a positive admissible
Green form for $\Pe(F_\C)$, (see \cite{BGS}, example 1. \@ 2.\@ (v)),
that is (\cite{BGS}. Prop.\@ 1.\@4.\@1)
\begin{equation} \label{Levgr}
d d^c [\Lambda_{\Pe(F)}] + \delta_{\Pe(F)} = \mu_M^p,
\end{equation}
and the harmonic projection of $\Lambda_{\Pe(E)}$ with respect to $\mu$ equals
\begin{equation} \label{Levproj}
H(\Lambda_{\Pe(E)}) = \sum_{n=1}^p \sum_{m=0}^{t-p} \frac{1}{m+n} \mu^{p-1} =
2(\sigma_t -\sigma_{p-1} - \sigma_{t-p}),
\end{equation}
that is
\begin{equation} \label{Levint}
\int_{\Pe(E_\C)} \Lambda_{\Pe(E)} \mu^{t-p+1} = 
\sum_{n=1}^p \sum_{m=0}^{t-p} \frac{1}{m+n}.
\end{equation}
(\cite{BGS},(1.\@4.\@1), (1.\@4.\@2))

\satz{Lemma} \label{normrelpr}
Denote by $|f|_\infty$ the sup norm of an element $f \in \check{E}_D$.
Then
\[ \log |f|_\infty - \frac{D}{2}\sum_{m=1}^t \frac{1}{m} \leq
   \int_{\Pe^t_\C} \log |f| \mu^t \leq \log |f|_{L^2} \leq \log |f|_\infty. \]
The first inequality is an equality iff $f$ is a power of a linear form.
\end{Satz}

\proof
\cite{BGS}, Prop.\@ 1.\@ 4.\@ 2, and formula (1.\@4.\@10).

\satz{Theorem} \label{Bez}
Let $\CX,\CY$ be effective cycles in $\Pe^t$ of codimension $p$, and $q$ 
respectively each being at most $t$, and assume that $p+q \leq t+1$,
and that $\CX$ and $\CY$ intersect properly. Then,
\[ h(\CX . \CY) \leq \deg (Y) h(\CX) + \deg(X) h(\CY) +
   (c_1-\sigma_t) \deg(\CX) \deg (\CY). \]
\end{Satz}

\proof
Follows immediately from Proposition \ref{vorbezout}.

%Let now $f \in \Gamma(\Pe^t,O(D))$. For any $\theta \in \Pe^t(\C)$ consider the
%restriction map
%\[ r_\theta: \Gamma(\Pe^t,O(D)) = \sym^D(\check{V}(\C)) 
%   \to \la \theta \ra, \]
%where $\la \theta \ra$ denotes the line in $\C^{t+1}$ corresponding to 
%$\theta$. If $F_D$ is the sub vector bundle of $\Pe^t(\sym^D(\check{D}))$
%over $\Pe^t$ whose fibre $F_{D,\theta}$ at $\theta$ is 
%$\mbox{Ker}(r_\theta)$, there is
%an exact sequence of vector bundles
%\[ 0 \longrightarrow F_D \longrightarrow \pi^*(\check{E}_D) 
%   \stackrel{r}{\longrightarrow} O(D) \longrightarrow 0. \]
%If $\theta$ is a representative of $\la \theta \ra$ of length one, then the 
%restriction of $f \in \check{E}_D$ to $\la \theta \ra$ 
%which is the same as the orthogonal projection of $f$ to the orthogonal
%complement of $F_{D,\theta}$
%has length $|f(\theta \otimes \cdots \otimes \theta)|$. This last expression 
%will also be simply denoted by $\la f | \theta \ra$.

\subsection{Flag varieties}

For $r \leq t$ and
$0 < q_1<q_2<\cdots < q_r  < q_{r+1}= t+1, q_i \in \N$, 
let $\CF=\CF_{(q_1,\ldots,q_r)}$ be the flag variety over $\spec \Z$ which 
assigns to each field $k$ the set of flags
\[ F: \quad 0 \subset V_1 \subset \cdots \subset V_r \subset k^{t+1}, \]
where each $V_i$ is a $q_i$-dimensional subspace of $k^{t+1}$;
denote by $d+1 = 1+\sum_{i=1}^r(q_{i+1}-q_i)q_i$ the dimension 
of $\CF$. 

We abbreviate $(q_1,\ldots,q_r)$ by $\bar{q}$.
On $\CF_{\bar{q}}$, there are the canonical quotient bundles
$Q_i, i=1,\ldots,r$, whose highest exterior powers $\CL_i, i=1,\ldots,r$
are ample generators of $\mbox{Pic}(\CF_{\bar{q}})$. Set
$\CL =\CL_{\CF}:= \otimes_{i=1}^r \CL_i$. A hermitian product on $\C^{t+1}$ 
canonically induces a hermitian metric on the base extensions
$L_i$ and $L$. We will
always assume that the canonical hermitian product on $\C^{t+1}$ was chosen, and
$U(t+1)$ the unitary group with respect to this product.
Set $\mu = \mu_F := c_1(\bar{L})$.

Since $U(t+1)$ operates transitively on $F_\C$, there is a unique 
$U(t+1)$-invariant volume form $\omega_{F}$ on $F_{\C}$ that is 
positive with respect
to the canonical orientation on $F_\C$ and gives $F_\C$ the
volume one. For an arbitrary point
$P_0 \in F(\C)$, by \cite{BGS}, 6.2, example (ii), there is a positive
green form of log type $g_0$ for $P_0$ such that 
$d d^c g_0 + \delta_{P_0} = \omega_F$.
For each point $P \in F(\C)$ choose an $h_p \in U(t+1)$ such that
$P = h_p P_0$, and define $g_P := (h_P^{-1})^* g_0$.
Since the metric on $\bar{L}$, and thereby $\mu_F$ is $U(t+1)$-invariant,
$g_P$ is positive for every $P$, and
\[ c_3(\bar{q}) := \frac12 \int_{F(\bar{q})} g_P \; \mu_{F(\bar{q})} \]
does not depend on $P$.

On $E_D = \Gamma(F_{\bar{q}},L^{\otimes D})$ the space of global sections of 
$L^{\otimes D}$, there are the norms
\[ |f|_\infty = \sup_{P \in F} |f_P|, \quad 
   |f|_m := \left(\int_F |f|^m \omega_F \right)^{\frac1m}, \quad
   |f|_0 := \exp \left(\int_F \log |f| \omega_F\right). \]

\satz{Proposition} \label{normrel}
The above norms fulfill the relations
\[ ||f||_0 \leq ||f||_m \leq ||f||_\infty \leq  \exp(c_3(\bar{q}) D)  ||f||_0. \]
\end{Satz}

\proof
The first two inequalities are valid for every probability space. For the
third inequality, let $f \in E_D$, and $P$ be a point in $G$. Then, 
$-\log  |f|^2$ is a Green form for $\di (f)$. Let $g_p$ be
the positive green form for $P$ on $F$ from above. 
By the commutativity of the star product,
\[ [-\log |f|^2] \delta_P + [g_P] D \mu_F = 
   [-\log |f|^2] \omega_F + g_P \delta_{\di(f)} \quad 
   \mbox{mod Im} \partial + \mbox{Im} \bar{\partial}. \]
Integrating over $F$ gives
\[ -\log |f_P|^2 + 2c_3(\bar{q}) D = -\log |f|_0^2 + \int_{\di(f)} g_p. \]
As $g_P$ is positive, $\int_{\di(f)} g_P \geq 0$, hence
\[ \log |f_P|^2 - 2c_3 D \leq  \log |f|_0^2 \]
for every $P \in F$ which implies
\[ \log |f|_\infty \leq c_3D + \log |f|_0. \]

\subsection{Grassmannians}

With the notations of the previous section, assume that
$\bar{q} = (q)$ is a single number. 
Then $\CF_{\bar{q}} =  \CG_{t+1,q}$ or $\CG_q$ or $\CG$ for short is the
Grassmannian that assigns to each field $k$ the set of flags
$0 \subset V_q \subset k^{t+1}$ consisting of only one space $V_q$ of
dimension $q$; in particular 
$\Pe^t = \CG_{t+1,q}$. The Picard group of $\CG$ is generated by the determinant 
$\CL = \CL_G$ of the canonical quotient bundle $\CQ$, hence
$c_1(\CL)$ is a generator for $CH^1(\CG) \cong \Z$.
The intersection product of every effective
cycle $X \in Z_p(G_\Q)$ with $c_1(L)^p$ is nonzero and defined as
the degree of $X$. 
Further, $\mu_G = c_1(\bar{L})$ is the K\"ahler form for the K\"ahler metric
induced by the canonical metric on $\C^{t+1}$ and the harmonic
forms are exactly
the forms that are invariant under $U(t+1)$. Hence, the product
of two harmonic forms is again harmonic, and by Proposition
\ref{schprodukt}, $CH^*(\bar{\CG})$ is a subring of 
$\widehat{CH}^*(\CG)$. In particular
$H^{0,0}(G(\C))$ are the constant functions, and $H^{d,d}(G(\C))$ are
the multiples of the volume form $\omega_G= \omega_{F_{(q)}}$.

With a fixed complete flag
\[ \{ 0 \} \subset V_1 \subset \cdots \subset V_t \subset V_{t+1} = k^{t+1} \]
and numbers $1 \leq i_1 < \cdots < i_q \leq t+1$ define the Schubert cell
$S^\circ_{(i_1,\ldots,i_q)}(k)$ in $G_q$ as the set of subspaces $W$ of
dimension $q$ such that
\[ \dim V_{i_l} \cap W = l \quad \forall \; l =1,\ldots,q. \]
Then the closure $S_{(i_1,\ldots,i_q)}$ of the Schubert cell is the set of
subspaces $W$ of dimension $q$ such that
 \[ \dim V_{i_l} \cap W \geq l \quad \forall \; l =1,\ldots,q. \]
The Schubert cells are algebraic subvarieties with
$\dim S_{(i_1,\ldots,i_q)}= \sum_{l=1}^q (i_l - l)$ and
\[ G_q = \bigcup_{(i_1,\ldots,i_q)}  S_{(i_1,\ldots,i_q)}, \]
the union being disjoint.
Further,
\[ Ch(\CG) = \bigoplus_{1\leq i_1<\cdots < i_q \leq t+1} \Z [S_{(i_1,\ldots,i_q)}]. \]
In particular, for $p\leq q$ the Schubert cell
$G_{V_p} = S_{(1,\ldots,p,t+2-(q-p),\ldots,t+1)}$ is the sub Grassmannian
consisting of the sub spaces $W$ that contain $V_p$, and for
$p \geq q$ the Schubert cell $G^{V_p} = S_{(p+1-q,\ldots,p)}$ is the sub
Grassmannian of the sub spaces $W$ that are contained in $V_p$. Further,
$S_{(1,\ldots,q)}$ represents the unique $0$-dimensional cycle class and
$S_{(t+1-q,t+3-q,\ldots,t+1)}$ represents the unique $q(t+1-q)-1$-dimensional
cycle class.

Denote $\sigma_{(i_1,\ldots,i_q)} = [S_{(i_1,\ldots,i_q)}] \in CH(G_q)$,
and the harmonic projection of $\delta_{S_{(i_1,\ldots,i_q)}}$ with
$\eta_{i_1,\ldots,i_q}$.

\satz{Proposition} \label{gras}
\begin{enumerate}
\item 
If $\sigma_{(i_1,\ldots,i_q)}$ and $\sigma_{(j_1,\ldots,j_q)}$ have complementary 
dimension, that is 
$\dim \sigma_{(i_1,\ldots,i_q)} + \dim \sigma_{(j_1,\ldots,j_q)} =
\sum_{l=1}^q i_l-l + \sum_{l=1}^q j_l-l = q(t+1-q)$, then
\[ \sigma_{(i_1,\ldots,i_q)} . \sigma_{(j_1,\ldots,j_q)} = \left\{
   \begin{array}{ll} \sigma_{(1,\ldots,q)} & if \quad
                                  i_l+j_{q+1-l} = t+2 \quad \forall l=1,\ldots,q\\
                            0 & \mbox{otherwise}. \end{array} \right. \]

\item
The harmonic projections 
$\omega_p: CH^p(G_q)_\C \to H^{p,p}(G_q), [X] \mapsto H(\delta_X)$
combine to a ring isomorphism $\omega: CH(G_q)_\C \cong H(G_q)$.
In particular $H^{1,1}(G_q)$ and $H^{(t+1-q)q,(t+1-q)q}(G_q)$ are one dimensional.

\item
The Arakelov Chow group $CH^p(\bar{G}_q)$ is isomorphic to the direct sum \\
$CH^p(G_q)_\C \oplus H^{p-1,p-1}(G_q)$ via the isomorphism
\[ CH^p(G_q)_\C \oplus H^{p-1,p-1}(G_q) \to CH^p(\bar{G}_q), \quad
   ([X], \eta) \mapsto s(X) + a(\eta). \]

\item With 
$d:=\sum_{l=1}^q j_l-l = q(t+1-q) - \left(\sum_{l=1}^q i_l-l\right)=
\dim S_{(i_1,\ldots,i_q)}$,
\[ \deg \sigma_{(i_1,\ldots,i_q)} = \frac{d!}{(i_1-1)! \cdots (i_q-1)!} 
   \prod_{l<k} (i_k-i_l). \]

\item
The restriction of $\eta_{i_1,\ldots,i_q}$ to $S_{j_1,\ldots,j_k}$ equals zero unless 
$i_l+j_{q+1-l} = t+1$ for every $l=1,\ldots,q$,
in which case it equals the unique volume form $\omega$ on
$= S_{i_1,\ldots,i_q}$ such that $S$ has volume one with respect to $\omega$. 
Consequently,
\[ \int_G \eta_{i_1,\ldots,i_q} \eta_{j_1,\ldots,j_q} = 
   \int_{S_{i_1,\ldots,i_q}} \eta_{(j_1,\ldots,j_q)} =
   \delta_{(i_1,\ldots,i_q).(t+1-j_1,\ldots,t+1-j_q)}, \]
and
\[ \int_{S_{i_1,\ldots,i_q}} \mu_G^d = \deg S_{i_1,\ldots,i_q}=
   \frac{d!}{(i_1-1)! \cdots (i_q-1)!} \prod_{l<k} (i_k-i_l). \]
In particular, $\mu_G^{q(t+1-q)} = \deg G_q \; \omega_G = 
\frac{1! 2! \cdots (q-1)! (q(t+1-q))!}{(t+2-q)! \cdots (t+1)!} \; \omega_G$
\end{enumerate}
\end{Satz}

\proof
1. \cite{Fu}, p.271.

\vspace{2mm}

2. \cite{Mai}, Th\'eore\`me 2.2.1.

\vspace{2mm}

3. Follows from part 2 and the Definition of the section $s$.

\vspace{2mm}

4. \cite{Fu}, Example 14.7.11.

\vspace{2mm}

5. The restriction of $\eta = \eta_{j_1,\ldots,j_q}$ to $S$ is a 
$U(i_1) \times U(i_2-i_1) \times \cdots \times U(i_q-i_{q-1})$-invariant
volume form on $S$. Hence $\eta|_S = a \omega$ with $a \in \R$.
Since there is a Green current $g$ for $S_{j_1,\ldots,j_q}$ with
\[ d d^c g + \delta_{S_{j_1,\ldots,j_q}} = \eta, \]
one gets
\[ \int_S \eta = (d d^c g) (\eta) + (\delta_S . \delta_{S_{j_1,\ldots,j_q}})(1) =
   g (d d^c \eta) + (\delta_s .  \delta_{S_{j_1,\ldots,j_q}})(1), \]
which because of $d d^c \eta = 0$ and part 1 equals $1$.

For the second equality, shorten $(i_1,\ldots,i_q)=I$, and 
$i_1+\cdots+i_q=|I|$. Because of part 2,
\[\mu^d = c_1(\bar{L})^d = H(c_1(L)^d)= 
   H(\sigma_{t+1-q,t+3-q,\ldots,t+1})^d =  H(\sigma_{t+1-q,t+3-q,\ldots,t+1}^d). \]
Hence, with $\sigma_{t+1-q,t+3-1,\ldots,t+1}^d = 
\sum_{|I|=(t+1)(t+1-q)-d-(q+1)q/2}  a_I\sigma_I$, 
\[ \int_S \mu^d = \sum_{|I|=(t+1)(t+1-q)-d-(q+1)q/2} \int_S a_I \eta_I, \] 
which by the above equals 
\[ a_{j_1,\ldots,j_q} = [S] . \; \sigma_{t+1-q,t+3-1,\ldots,t+1}^d =
  [S] . c_1(L)^d = \deg S, \]
and the claim follows from part 4.

\vspace{2mm}

For $\bar{q}=(q_1,\ldots,q_r)$ the canonical maps
$\varphi_{\bar{q},q_i}: \CF_q \to \CG_{q_i}=\CG_{t+1,q_i},,i=1,\ldots,r$ that 
forget every subspace except the $i$th, are projective bundle maps,
hence flat, proper, smooth, surjective, and projective.
For $\CF= \CG_{(1,\ldots,t)}$ the full flag variety, and
$g_p : \CF \to \CG_p, p =1,\ldots,t$ the form
$(g_p)_* \omega_F$ is $U(t+1)$ invariant, and
$\int_{G_p} (g_p)_* \omega = \int_F \omega = 1$, hence
$(g_p)_* \omega_F = \omega_p:= \omega_{G_p}$.

For $\bar{q} = (q,p)$, there is the correspondence

\special{em:linewidth 0.4pt}
\linethickness{0.4pt} 
\vspace{4mm}
\hspace{5cm}
\begin{picture}(160.00,40.00)
\put(68.00,28.00){\vector(-3,-2){24.00}}
\put(80.00,35.00){\makebox(0,0)[cc]{$\CF_{q,p}$}}
\put(30.00,5.00){\makebox(0,0)[cc]{$\CG_q$}}
\put(130.00,5.00){\makebox(0,0)[cc]{$\CG_p$}}
\put(54,28){\makebox(0,0)[cc]{$\varphi_q$}}
\put(106,28){\makebox(0,0)[cc]{$\varphi_p$}}
\put(92.00,28.00){\vector(3,-2){24.00}}
\end{picture}
\vspace{-11mm}
\begin{equation} \label{cor1} \end{equation}
\vspace{-2mm}

between $\CG_q=\CG_{t+1,q}$, and $\CG_p = \CG_{t+1,p}$. 
For later applications we will need the fact that the maps
$\varphi_q$ and $\varphi_p$ preserve the dimension of certain Schubert 
cells.

\satz{Lemma} \label{dimension}
With the above notations
\begin{enumerate}
\item
let $q=1$. Then, 
\[ \dim (\varphi_p (\varphi_1^{-1} S_{t+1-p})) = 
   \dim (\varphi_p(\varphi^{-1} \Pe^{t-p})) = 
   \dim (\varphi^{-1} \Pe^{t-p}), \]
and with $W$ a space of dimension $p$,
\[ \dim (\varphi_1 (\varphi^{-1}_p S_{1,\ldots,p}) =
   \dim (\varphi_1 (\varphi^{-1}_p W)) = \dim (\varphi^{-1} W). \]
\item
The maps $(\varphi_p)_* \circ \varphi_q^*,(\varphi_q)_* \circ \varphi_P^*$ map
harmonic forms to harmonic forms.
\end{enumerate}
\end{Satz}

\proof
1. For $[v] \in \Pe^t$ a point, by Proposition \ref{gras},
\[ \dim \varphi_1^{-1} ([v]) = \dim S_{(1,t+3-p, \ldots,t+1)} = (t+1-p)(p-1), \]
hence
\[ \dim \varphi_1^{-1} \Pe^{t-p} = \dim \Pe^{t-p} + (p-1) (t+1-p) = 
   t-p + (p-1)(t+1-p), \]
\[ \dim (\varphi_p (\varphi_1^{-1} \Pe^{t-p}) = \]
\[ \dim S_{t+1-p, t+3-p,\ldots t+1} = t+1-p-1 + \sum_{l=2}^p t+1+l-p -l =
   t-p + (p-1)(t+1-p), \]
and also,
\[ \dim (\varphi^{-1}_p W) = \dim S_{(p)} =1\cdot(p-1), \]
and
\[ \dim (\varphi_1 (\varphi^{-1}_p W)) = \dim \Pe(W) = p-1. \]

\vspace{2mm}

2. Follows from the facts, that the harmonic forms on $G_p,G_q$ are exactly
the forms that are invariant under $U(t+1)$ and that $\varphi_p,\varphi_q$
are $U(t+1)$-invariant.

\satz{Corollary} \label{VX}
Let $X$ be an effective cycle of pure codimension $t+2-p$ in 
$\Pe^t$.
Then, $(\varphi_p)_* \varphi_1^* X$ is an effective cycle that has pure 
codimension $1$ and the same dimension as $\varphi_1^*X$. We will
use the abbreviation $V_X = (\varphi_p)_* \varphi_1^* X$.
\end{Satz}

\proof
Follows from the Lemma together with $[X] = \deg X [S_{t+1-p}]$.

\vspace{2mm}

With $\Pe(W) \subset \Pe^t(\C)$ a subspace of dimension $q$, define the map
\[ \pi_{W^\bot}: \Pe^t \setminus W \to W^\bot, \quad
   [v+w] \mapsto [w], \quad v \in \Pe(W), w \in \Pe(W^\bot). \]
If $X$ is an effective cycle in $\Pe^t$ whose support does not meet
$\Pe(W)$, define 
\begin{equation} \label{XW}
X_W := \overline{\pi^*(\pi_* X)}.
\end{equation}

\satz{Proposition} \label{doppelcor}
\begin{enumerate}
\item
For $p>q$, consider the 3 correspondences on complex manifolds

\special{em:linewidth 0.4pt}
\linethickness{0.4pt} 
\vspace{4mm}
\hspace{1cm}
\begin{picture}(200.00,40.00)
\put(38.00,28.00){\vector(-3,-2){24.00}}
\put(50.00,35.00){\makebox(0,0)[cc]{$F_{1,p}$}}
\put(0.00,5.00){\makebox(0,0)[cc]{$\Pe^t$}}
\put(100.00,5.00){\makebox(0,0)[cc]{$G_p$}}
\put(24,28){\makebox(0,0)[cc]{$\varphi_1$}}
\put(76,28){\makebox(0,0)[cc]{$\varphi_p$}}
\put(62.00,28.00){\vector(3,-2){24.00}}

\put(168.00,28.00){\vector(-3,-2){24.00}}
\put(180.00,35.00){\makebox(0,0)[cc]{$F_{1,p-q}$}}
\put(130.00,5.00){\makebox(0,0)[cc]{$\Pe^t$}}
\put(230.00,5.00){\makebox(0,0)[cc]{$G_{p-q}$}}
\put(154,28){\makebox(0,0)[cc]{$\psi_1$}}
\put(211,28){\makebox(0,0)[cc]{$\psi_{p-q}$}}
\put(192.00,28.00){\vector(3,-2){24.00}}

\put(303 ,28.00){\vector(-3,-2){24.00}}
\put(315.00,35.00){\makebox(0,0)[cc]{$F_{p-q,p}$}}
\put(265.00,5.00){\makebox(0,0)[cc]{$G_{p-q}$}}
\put(365.00,5.00){\makebox(0,0)[cc]{$G_p$}}
\put(284,28){\makebox(0,0)[cc]{$\bar{\varphi}_{p-q}$}}
\put(341,28){\makebox(0,0)[cc]{$\bar{\varphi}_p$,}}
\put(327.00,28.00){\vector(3,-2){24.00}}

\end{picture}
\vspace{-7mm}
\[ \]

let $\Pe(W) \subset \Pe^t$ be a $q$-dimensional subspace,
$X \in Z_{eff}^p(\Pe^t)$ an effective cycle whose support does not
intersect $\Pe(W)$, and $G_W \subset G_p$  the Schubert cell consisting
of the subspaces that contain $W$. Then
\begin{equation} \label{GWgl}
(\bar{\varphi}_{p-q})_* \bar{\varphi}_p^*
(G_W . (\varphi_p)_* (\varphi_1^* X)) = (\psi_{p-q})_* (\psi_1^* X_W),
\end{equation}
\[ \dim (\bar{\varphi}_{p-q})_* \bar{\varphi}_p^*
   (G_W . (\varphi_p)_* (\varphi_1^* X)) =
   \bar{\varphi}_p^*
   (G_W . (\varphi_p)_* (\varphi_1^* X)). \]

\item 
For $p \leq q$, with the correspondences on complex manifolds

\special{em:linewidth 0.4pt}
\linethickness{0.4pt} 
\vspace{4mm}
\hspace{1cm}
\begin{picture}(200.00,40.00)
\put(38.00,28.00){\vector(-3,-2){24.00}}
\put(50.00,35.00){\makebox(0,0)[cc]{$F_{1,p}$}}
\put(0.00,5.00){\makebox(0,0)[cc]{$\Pe^t$}}
\put(100.00,5.00){\makebox(0,0)[cc]{$G_p$}}
\put(24,28){\makebox(0,0)[cc]{$\varphi_1$}}
\put(76,28){\makebox(0,0)[cc]{$\varphi_p$}}
\put(62.00,28.00){\vector(3,-2){24.00}}

\put(168.00,28.00){\vector(-3,-2){24.00}}
\put(180.00,35.00){\makebox(0,0)[cc]{$F_{1,p+q}$}}
\put(130.00,5.00){\makebox(0,0)[cc]{$\Pe^t$}}
\put(230.00,5.00){\makebox(0,0)[cc]{$G_{p+q}$}}
\put(154,28){\makebox(0,0)[cc]{$\psi_1$}}
\put(211,28){\makebox(0,0)[cc]{$\psi_{p+q}$}}
\put(192.00,28.00){\vector(3,-2){24.00}}

\put(303,28.00){\vector(-3,-2){24.00}}
\put(315.00,35.00){\makebox(0,0)[cc]{$F_{p+q,p}$}}
\put(265.00,5.00){\makebox(0,0)[cc]{$G_{p+q}$}}
\put(365.00,5.00){\makebox(0,0)[cc]{$G_p$}}
\put(284,28){\makebox(0,0)[cc]{$\bar{\varphi}_{p+q}$}}
\put(341,28){\makebox(0,0)[cc]{$\bar{\varphi}_p$,}}
\put(327.00,28.00){\vector(3,-2){24.00}}

\end{picture}
\vspace{-7mm}
\[ \]
let $\Pe(F) \subset \Pe^t$ be a $q$-codimensional subspace,
$X \in Z_{eff}^p(\Pe^t)$ an effective cycle intersecting $\Pe(F)$ properly,
and $G^F \subset G_p$ the sub Grassmannian of subspaces that are contained
in $F$. Then
\begin{equation} \label{GFgl}
(\bar{\varphi}_{p+q})_* \bar{\varphi}_p^* 
(G^F . (\varphi_p)_* (\varphi_1^* X)) = (\psi_{p+q})_* (\psi_1^* (X.\Pe(F)),
\end{equation}
where $G^F \subset G_p$ is the subset of spaces being contained in $F$.
Further,
\[ \dim (\bar{\varphi}_{p+q})_* \bar{\varphi}_p^* 
   (G^F . (\varphi_p)_* (\varphi_1^* X)) =
   \dim \bar{\varphi}_p^* 
   (G^F . (\varphi_p)_* (\varphi_1^* X)). \]
\end{enumerate}
\end{Satz}

\proof
For $X$ a variety, the equalities of proof for the equalities (\ref{GWgl}) and
(\ref{GFgl}) are easy exercises in linear algebra, thus hold for arbitrary
cycles by linearity. For the equalities of dimensions,

1. Since $G_W$ is isomorphic to the Grassmannian of $p-q$-dimensional subspaces
in $\C^{t+1-q}$ and 
$\mbox{codim} V_X = \mbox{codim} (\varphi_p)_*(\varphi_1^*X)) =1$, we have
$\dim G_W . (\varphi_p)_*(\varphi_1^*X)) = (p-q)(t+1-p)-1$. Further, the
fibre of $F_{p-q,p}$ above each point $E$ in $G_p$ consists of the $(p-q)$-
dimensional subspaces of $E$, and has thus dimension $q(p-q)$. 
Hence
\[ \dim \bar{\varphi}^*_p(\dim G_W . (\varphi_p)_*(\varphi_1^*X)) = \]
\[ (p-q)(t+1-p)-1 + \dim G_{p-q} = (p-q)(t+1-p)-1 + q(p-q) =
   (p-q)(t+1+q-p)-1. \]
On the other hand, 
\[ \dim V_{X_W}=\dim (\psi_{p-q})_*(\psi_1^*X_W) = \dim G_{p-q}-1 =
   (t+1+q-p)(p-q)-1, \]
hence, the two dimensions are equal.

\vspace{2mm}

2. Again, $\dim (G^F .(\varphi_p)_* (\varphi_1^* X)) = \dim G^F-1 =
p(t+1-p-q)-1$,
and a fibre of $\bar{\varphi}_p$ has dimension $q(t+1-p-q)$. Hence,
$\dim \bar{\varphi}_p^* (G^F . (\varphi_p)_* (\varphi_1^* X)) =
(p+q)(t+1-p-q)-1$.
Since also
\[ \dim (\psi_{p+q})_* (\psi_1^* (X.\Pe(F)) = \dim V_{X.\Pe(F)} = G_{p+q}-1 =
   (p+q)(t+1-p-q)-1, \]
the dimensions are equal.

\subsection{Chow divisor}

Let $\check{\Pe}^t$ be the dual projective space.
The $p$ projections 
$pr_i : (\check{\Pe^t})^p \to \check{\Pe^t}$, define line bundles
$O_i(1) = pr^*_i(O(1))$, and
$O(D_1, \ldots, D_p) = O_1(1)^{\otimes D_1} \otimes \cdots \otimes 
O_p(1)^{\otimes D_p}$.

A dual inner product on $\check{\Pe^t}$ defines metrics on $O(D)$, and
$O(D_1, \ldots D_p)$, and a K\"ahler metric on $\check{\Pe^t}$, and
$(\check{\Pe}^t)^p$. The corresponding K\"ahler forms are
$\check{\mu} = c_1(\overline{O(1)})$, and
$\bar{\mu} = \check{\mu}_1 + \cdots + \check{\mu}_p =
pr_1^* \check{\mu} + \cdots + pr_p^* \check{\mu} =
c_1(\overline{O(1,\ldots,1)})$. The space of harmonic forms on
$(\check{\Pe}^t)^p$ is the linear span of the products
$\prod_{i=1}^p \check{\mu}_i^{k_i}, \; 1 \leq k_i \leq t \; \forall i=1,\ldots,p$.

Let $\delta: \check{\Pe}^t \to (\check{\Pe}^t)^p$ be the diagonal, and
define the correspondence

\special{em:linewidth 0.4pt}
\linethickness{0.4pt} 
\vspace{4mm}
\hspace{5cm}
\begin{picture}(160.00,40.00)
\put(68.00,28.00){\vector(-3,-2){24.00}}
\put(80.00,35.00){\makebox(0,0)[cc]{$\CC$}}
\put(30.00,5.00){\makebox(0,0)[cc]{$(\Pe^t)^p$}}
\put(130.00,5.00){\makebox(0,0)[cc]{$(\check{\Pe}^t)^p$}}
\put(54,28){\makebox(0,0)[cc]{$f$}}
\put(106,28){\makebox(0,0)[cc]{$g$}}
\put(92.00,28.00){\vector(3,-2){24.00}}
\end{picture}
\vspace{-11mm}
\begin{equation} \label{chow} \end{equation}
\vspace{-2mm}

where $\CC$ is the sub scheme of $(\Pe^t)^p \times (\check{\Pe^t})^p$ 
assigning to each $t+1$ dimensional vector space $V$ over a field $k$ the set  
\[ \{ (v_1, \ldots, v_p,\check{v}_1,\ldots \check{v}_p | v_i \in V,\;
   \check{v}_i \in \check{V}, \; \check{v}_i(v_i) = 0, \; 
   \forall i = 1,\ldots p. \} \] 
The maps
$f:\CC \to (\Pe^t)^p, g:\CC \to (\check{\Pe}^t)^p$ are just the restrictions
of the projections. They are flat, projective, surjective,
and smooth.

Let $\CX \in Z^{t+1-p}_{eff}(\Pe^t)$,
and define the Chow divisor $Ch(\CX) \subset (\check{\Pe}^t)^p$ as
$Ch(\CX) := g_* \circ f^* \circ \delta_* (\CX)$.

\satz{Proposition} \label{chschieb}
\begin{enumerate}

\item
The Chow divisor has codimension one; it is the divisor corresponding
to a global section 
$f_X \in \Gamma((\check{\Pe}^t)^p,O(\deg X,\ldots,\deg X))$ such that
\[ d d^c [- \log |f_X|^2] + \delta_{Ch(X)} = \deg X \bar{\mu}. \]
Consequently, $-\log |f_X|^2$ is an admissible Green form of log type
for $Ch(X)$, and for all $i=1, \ldots, p$ 
the multi degrees of $Ch(X)$, that is the numbers
\[ c_1(O_1(1))^t . \cdots . c_1(O_{i-1}(1))^t . c_1(O_i(1))^{t-1} \cdot
   c_1(O_{i+1}(1))^t . \cdots . c_1(O_p(1))^t . [Ch(X)] \]
all equal $\deg X$, i.\@ e.\@ 
$[Ch(X)] = (\deg X, \ldots, \deg X) \in \Z^p = CH^1((\check{\Pe}^t)^p)$.

Further $\dim X = \dim \delta(X)$, and 
$\dim g^* \delta_* X = \dim f(g^* \delta_* X)$.

\item
If $\Pe(W)$ does not meet $X$, then $\Pe(\check{W})$, and $Ch(X)$
intersect properly.
Further if $\check{w}_1, \ldots, \check{w}_p \in \check{W}$ are $p$ linearly 
independent vectors, and $\Pe(W)$ the intersection of their kernels, then
$\dim g^{-1} (\check{w}_1, \ldots \check{w}_p) = 
\dim f(g^{-1} (\check{w}_1, \ldots \check{w}_p))$, and
$\delta^* f_* g^* (\check{w}_1, \ldots \check{w}_p) = \Pe(W)$.

\item
The maps $g_* \circ f^*,f_* \circ g^*$ map harmonic forms to harmonic 
forms.

\item
For $k \leq t+1-p$ let $\check{v}_1, \ldots,\check{v}_k$ be orthonormal
vectors, $E_i= \mbox{ker} (\check{v}_i), i=1,\ldots,k$, and
$E = E_1 \cap \cdots \cap E_k$. Let further $g$ be a normalized green
form for $E_1 \times \cdots \times E_k \times (\Pe^t)^{t+1-p-k}$ in
$(\Pe^t)^p$, and 
$i^{-1}(g)$ its restriction to $\Pe^t$ via the diagonal embedding $\delta$.
The number
\[ c_5 := \int_{\Pe^t} i^{-1}(g) \mu^{t+1-k} \]
only depends on $t,p,k$, and $i^{-1}(g) - c_5 \mu^t$ is a normalized
green form for $E$.

\end{enumerate}
\end{Satz}

\proof
1. One only has to check that a generator $[\Pe(V)] \in CH^{t+1-p}(\Pe^t)$
by $g_* f^* \delta_*$ is mapped to $(1, \ldots,1) \in CH^1((\check{\Pe}^t)^p)$,
which is obvious.

\vspace{2mm}

2. is obvious.

\vspace{2mm}

3. Follows from the fact that the harmonic forms on $(\Pe^t)^p$ and
$(\check{\Pe})^p$ are exactly the forms invariant under $U(t+1)^p$, and
the $U(t+1)$ equivariance of $f$ and $g$.

\vspace{2mm}

4. Follows immediately from the fact that $U(t+1)$ acts transitively
on the set orthonormal $k$-tupels.

\vspace{2mm}

The canonical quotient bundle $Q$ on $(\check{\Pe^t})$ carries a canonical
metric as well. Let $\hat{c}(\bar{Q})$ be its total arithmetic Chern class.
(See \cite{SABK}). Define the height of a divisor $D$ in $(\check{\Pe^t})^p$ as
\[ h(D) = \pi_*(\hat{c}(\prod_{i=1}^p\bar{Q}_i)|D). \]

\satz{Proposition}
For all effective cycles $X \in Z^{t+1-p}_{eff}(\Pe^t)$,
\[ h(\CX) = h(Ch(\CX)). \]
\end{Satz}

\proof
\cite{BGS}, Theorem 4.\@3.\@2.

\vspace{2mm}

With $(d_1, \ldots d_p) \in \N^p$ the space
$E_{d_1,\ldots ,d_p} = \Gamma(\check{\Pe}^t,O(d_1,\ldots,d_p))$ carries
norms \\ $||\cdot||_0, ||\cdot||_r, ||\cdot||\infty, r \in \R^{>0}$ just
like $\Gamma(\Pe^t,O(D))$, and the analogous estimates hold.

\satz{Lemma} \label{chnab}
There is a positive constant $c_7$, only depending on $t$ such that
for any $f \in E_{d_1,\ldots,d_p}, r \in \R^{>0}$,
\[ \log ||f||_\infty - c_7 \sum_{i=1}^p d_i \leq \log ||f||_0 \leq 
   \log ||f||_r \leq \log ||f||_\infty. \] 
\end{Satz}

\proof
\cite{BGS}, Corollary 1.\@ 4.\@ 3.

\vspace{2mm}

%\vspace{2mm}

%Let now $\Pe(W),\Pe(V)$ be arbitrary subspaces of of dimension $r \leq q$,
%$s \geq q$ of $\Pe^t$, and consider the subvariety $G_W,V$ of $G$
%consisting of the $q$ dimensional subspaces of $\Pe^t$ containing $\Pe(W)$
%and being contained in $\Pe(V)$. This variety is canonically isomorphic
%to $G_{s,q-r}$ the set of $q-r$-dimensional subspaces of $V$.

%\satz{Corollary} \label{abstrel}
%With notations as above 
%\[ -\frac12 \int_G \int_{\Pe(F)} g_x \mu_G^d \geq 
%   -\frac12 \stackrel{sup}{F \in G} \int_{\Pe(F)} g_X \geq 
%   -\frac12 c(p,t) D -\frac12 \int_G \int_{\Pe(F)} g_x \mu_G^d. \]
%\end{Satz}
%
%\proof
%With $f_{f_X}$ as in Proposition \ref{glschnitt}, this follows from
%\[ -\frac12 \int_G \int_{\Pe(F)} g_x \mu_G^d = 
%    \int_G \log |f_{V_X}| +c, \]
%\[ -\frac12 \stackrel{sup}{F \in G} \int_{\Pe(F)} g_X = 
%   \stackrel{sup}{F \in G} \log |f_{V_X}(\Pe(F))|+c, \]
%and the previous proposition.

%\satz{Proposition}
%If $p <q$, then there is are positive constants 
%$c_0,c_3,c_\infty$ only depending on $p,q$, and $t$ such that for
%every $f \in E_D$,
%\[ \log ||f||_0 \geq -c_0 D, \quad \log ||f||_2 \geq c_3 D, \quad
%   \log ||f||_\infty \geq -c_\infty D. \]
%\end{Satz}

%\proof
%Because of the previous Proposition, the claim has to be proved
%only for one of the norms.

\section{The algebraic distance}

In this section all varieties, and cycles are assumed to be 
defined over $\C$. All integrals over
sub varieties of smooth projective varieties are defined via
resolutions of singularities (cp.\@ \cite{SABK}, II.\@1.\@2.)

\subsection{Definitions and fundamental properties}

%The main ingredient of the proof of the metric B\'ezout Theorem is that
%on Arakelov varieties, there is a concept of algebraic of effective
%cycles that measures the difference between the height of their
%intersection product and their heights and degrees.
Let $X$ be a smooth projective variety over $\C$,
and fix a K\"ahler metric with K\"ahler form $\mu$ on $X$.
For $p+q \leq t+1$, define the pairing
\[ Z^p(X_\C) \times Z^q(X_\C) \to \CD^{t,t}(X_\C), \quad 
   (Y,Z) \mapsto (Y|Z) := [g_Y] (\delta_Z - H(\delta_Z)) \mu^{t+1-p-q} \]
on the cycle group, where $g_Y$ is an admissible Green form of log type
for $Y$.

\satz{Lemma and Definition}
The above pairing is well defined and symmetric modulo 
$\mbox{Im} \partial + \mbox{Im} \bar{\partial}$. 
If $X$ and $Y$ intersect properly, the algebraic distance
\[ D(X,Y) := -\frac12 \int_X  (Y|Z) \]
is finite.  
\end{Satz}

\proof
Let $g_Y'$ be another admissible Green form for $Y$. By  Proposition 
\ref{sequ},
$g_Y - g_Y'$ modulo $\mbox{Im} \partial + \mbox{Im} \bar{\partial}$
equals a harmonic form $\eta$. Thus,
\[ [g_Y] (\delta_Z - H(\delta_Z)) \mu^{t+1-p-q} -
   [g'_Y] (\delta_Z - H(\delta_Z)) \mu^{t+1-p-q} = \]
\[ (\delta_Z - H(\delta_Z)) \eta \mu^{t+1-p-q} =
   - (d d^c [g_Z]) \eta \mu^{t+1-p-q} = \]
\[ d^c [g_Z] d (\eta \mu^{t+1-p-q}) 
   \quad \mbox{mod Im} d \subset \mbox{Im} 
   \partial + \mbox{Im} \bar{\partial}. \]
The last expression equals zero, since harmonic forms are contained 
in the kernel of $d$. It follows that the pairing is well defined.

For the symmetry, we have to prove
\[ [g_Y] (\delta_Z - H(\delta_Z)) \mu^{t+1-p-q} =
   [g_Z] (\delta_Y - H(\delta_Y)) \mu^{t+1-p-q}
   \quad \mbox{mod Im} \partial + \mbox{Im} \bar{\partial} \]
for admissible Green forms $g_Y,g_Z$.
This is equivalent to
\[ [g_Y] \delta_Z + [g_Z] H(\delta_Y) = [g_Z] \delta_Y + [g_Y] H(\delta_Z) 
   \quad \mbox{mod Im} \partial + \mbox{Im} \bar{\partial}, \]
which just is the commutativity of the star product of Green currents.

The last claim of the Lemma follows from the fact that $g_Y$ is of
log type along $Y$.

\vspace{2mm}

One immediately observes

\satz{Fact}
If one of the cycles $X,Y$ is the zero cycle, then $D(X,Y) = 0$.
For $p+q \leq t+1$ the map $D: Z^p_{eff}(X) \times Z^q_{eff}(X) \to \R$ is
bilinear on the subset on which it is defined, i.\@e.\@ for properly
intersecting cycles.
\end{Satz}

\satz{Remark} \label{remark}
If $g_Y$ is normalized, then
\[ D(Y,Z) = -\frac12 \int_{X(\C)} g_Y \delta_Z. \]
If $X$ is the base extension of an arithmetic variety $\CX$ and
on $X(\C)$ the product of harmonic forms is again harmonic, then 
$H(g_Y (\delta_Z - H(\delta_Z))) = H(g_Y \delta_Z)$, and by Proposition
\ref{hsch}.\@2,
\[ h(\CY . \CZ) = 
   \pi_*(\hat{c}_1(\CL)^{t-p-q+1} s(\CY) . s(\CZ)) + D(Y,Z). \] 
In case $X = \Pe^t$, Proposition \ref{vorbezout}.1 reformulates as
\[ h(\CY . \CZ) = \deg Y h(\CZ) + \deg Z h(\CY) + D(Y,Z)
   - \sigma_t \deg Y \deg Z, \]
and Theorem \ref{Bez} reformulates as
\[ D(X,Y) \leq \bar{c}_1 \deg X \deg Y. \]
\end{Satz}

\satz{Lemma} \label{Gnoben}
Let $Y \in Z_{eff}^1(G)$, with $G=G_{t,p}$ the Grassmannian, and 
$V \in G(\C)$ a point not contained in the support of $Y$. Then, with
$V \subset \C^{t+1}$ a subspace of dimension $p$, i.\@ e.\@ a point in $G$,
\[ D(Y,V) \leq c_3(t,p) \deg_{L_G} Y. \]
\end{Satz}

\proof
Let $f \in \Gamma(G,L_G)$ such that $Y = \di f$. Since with $d= \dim G$,
\[ D(Y,V) = \log |f_V| - \int_G \log |f| \mu_G^d =
            \log |f_V| - \log ||f||_0, \]
the Lemma immediately follows from Proposition \ref{normrel}.

\vspace{2mm}

The following two propositions supply the essential techniques for calculating
algebraic distances.

Let $Z,W$ be projective K\"ahler varieties of dimensions $r,s$ 
with K\"ahler structures $\mu_Z,\mu_W$, and $C$ a projective regular
variety over $\C$ of pure dimension $d$, and consider a correspondence

\special{em:linewidth 0.4pt}
\linethickness{0.4pt} 
\vspace{4mm}
\hspace{5cm}
\begin{picture}(160.00,40.00)
\put(68.00,28.00){\vector(-3,-2){24.00}}
\put(80.00,35.00){\makebox(0,0)[cc]{$C$}}
\put(30.00,5.00){\makebox(0,0)[cc]{$Z$}}
\put(130.00,5.00){\makebox(0,0)[cc]{$W$}}
\put(54,28){\makebox(0,0)[cc]{$f$}}
\put(106,28){\makebox(0,0)[cc]{$g$}}
\put(92.00,28.00){\vector(3,-2){24.00}}
\end{picture}
\vspace{-11mm}
\begin{equation} \label{cor2} \end{equation}
\vspace{-2mm}

with $f,g$ flat, surjective and projective. 
For $X \in Z_{eff}(Z), Y \in Z_{eff}(W)$
define $C_* (X) := g_*(f^*(X)), C^*(Y) := f_*(g^*(Y))$.

%\[ X' \longleftarrow X \longrightarrow X'' \]
%is a correspondence, with $f: X \to X'$, and $g: X \to X''$ flat, and dominant
%and mapping harmonic forms to harmonic forms, then for
%effective cycles $Y$ in $X'$, and $Z$ in $X''$ of appropriate dimension, 

\satz{Proposition (Functoriality)} \label{Dfunk1}
In the above situation, let $q \geq p$ and
$X \in Z_{eff}^p(Z), Y \in Z_{eff}^q(W)$ be such that
$C_*(X)$ and $C^*(Y)$ are both nonzero, and the intersections of
$X$ with $C^*(Y)$ and of $C_*(X)$ with $Y$ are both proper.

\begin{enumerate}
\item
If $C^* = f_* \circ g^*$ maps harmonic forms to harmonic forms,
then $C_* = g_* \circ f^*$ maps normalized Green forms to normalized
Green forms.

\item
If $p+q=d+1$ and $C^* = f_* \circ g^*$ maps harmonic forms to harmonic
forms, then
\[ D(X,C^*Y) = D(C_*X,Y). \]

%\item
%If $p+q < d+1$, $C^* = f_* \circ g^*$ maps harmonic forms to harmonic forms
%and additionally $C_* \mu_X = \mu_Y$, still
%\[ D(X,C^*Y) = D(C_*X,Y). \]
\end{enumerate}
\end{Satz}

\proof
1. 
Let $g_X$ be a normalized Green form for $X \in Z_{eff}^p(Z)$, and \\
$\eta \in H^{d+1-p,d+1-p}(W)$. Since by
\cite{SABK}, Lemma II.\@2 (ii), $g_*[f^*g_Z] = [g_* (f^*g_Z)]$,
\[ \int_W C_* g_X \; \eta = \int_W g_* (f^* g_X) \; \eta = 
   \int_C f^*g_X g^* \eta = \int_Z g_X f_*(g^* \eta)= 0, \]
as $f_*\circ g^* \eta$ is harmonic, that is $C_* g_X$ is orthogonal to
$H^{d+1-p,d+1-p}(W)$.

\vspace{1mm}

2. By Remark \ref{remark}, with $g_X$ a normalized Green form for $X$,
\[ D(X,C^* Y) = -\frac12\int_{C^* Y} g_X =
   [\C(g^{-1}Y):\C(f(g^{+1}(Y)))] \int_{f(g^{-1}(Y))} g_X = \]
\[ \int_{g^{-1}(Y)} f^* g_X, \]
which by Fubini's Theorem equals
\[ \int_Y (g_* \circ f^*) g_X. \]
As, by the first part $(g_* \circ f^*) g_z$ is a normalized Green 
form, this in turn equals $D(C_* X ,Y)$.

%\vspace{2mm}
%3.  a normalized Green form $g_Z$ for $Z$,
%\[ D(f_* Y,Z) = \int_{f_* Y} g_Z \mu_{X'}^{q-p} =
%   [\C(Y):\C(f(Y))] \int_{f(Y)} g_Z \mu_{X'}^{q-p} = 
%   \int_Y f^* (g_Z \mu_{X'}^{q-p})= \]
%\[ \int_Y f^* g_Z \mu_X^{q-p}. \]
%Again by the Lemma, this equals $D(Y,f^*Z)$. 

%\vspace{3mm}

%2. By the Ramark,
%\[ D(f^*Y,Z) = \int_{f^*Y} g_Z. \]
%By the definition of $f_* g_Y$, this equals
%\[ \int_Z f_* [g_Y] = 
%   \left(\int_Z f_* [g_Y] - \int_X f_* g_Y \; \omega_Z \right) + 
%   \int_X f_* [g_Y] \; \omega_Z, \]
%which by definition equals
%\[ D(Y,f^*Z) = D(f_*Y,Z) + \int_{X'} f_* g_Y \; \omega_Z. \]

\vspace{2mm}

\satz{Proposition} \label{vergl}
Let $W$ be a projective algebraic K\"ahler variety of dimension $t$
with K\"ahler form $\mu$, and for
$p+q+r \leq t+1$, let $X,Y,Z$  be effective cycles of pure codimension 
$p,q,r$  on $W$
such that $X . Y, Y . Z, X . Y, Y . Z . W$ are of pure 
codimension $p+q,q+r,p+r,p+q+r$ respectively. Then,

\begin{enumerate}

\item
If $g_Y,g_Z$ are admissible Green forms for $Y$ and $Z$
\[ [g_Y] \wedge \delta_{X . Z} + 
   H(\delta_Y) \wedge [g_Z] \wedge \delta_X=
   \delta_{X . Y} \wedge [g_Z] + H(\delta_Z) \wedge [g_Y] \wedge\delta_X \quad
   \mbox{mod Im} \partial + \mbox{Im} \bar{\partial}. \]
In particular, for $W= \Pe^t$ projective space,
\[ [g_Y] \wedge \delta_{X . Z} + 
   \deg Y \mu^q \wedge [g_Z] \wedge \delta_X=
   \delta_{X . Y} \wedge [g_Z] + \deg Z \mu^r \wedge [g_Y]\wedge\delta_X 
   \quad\!\!\!
   \mbox{mod Im} \partial + \mbox{Im} \bar{\partial}. \]

\item
\[ D(Y, X . Z) -\frac12 \int_X [g_Z] H(\delta_Y) \mu^{t+1-p-q-r} = \]
\[ D(X . Y,Z) - \frac12 \int_X [g_Y] H(\delta_Z) \mu^{t+1-p-q-r}. \]
In particular, in projective space,
\[ D(Y,X . Z) + \deg Y D(X,Z) =
   D(X . Y,Z) + \deg Z D(X,Y). \]
\end{enumerate}
\end{Satz}

\proof
1.
\cite{GS1}, Theorem 2.\@2.\@2.

\vspace{2mm}

2. In 1, let $g_Y,g_Z$ be the $\mu$-normalized Green forms of $Y$, and $Z$. 
Multiplying the equality of 1 
with $\mu^{t+1-p-q-r}$, integrating over $X$, and dividing
by $-2$ leads the first equality. The second equality follows from
$H(\delta_Y) = \deg Y \mu^{t+1-q}, H(\delta_Z) = \deg Z \mu^{t+1-r}$ in
case of projective space, and $g_Y,g_Z$ normalized and hence admissible.

\vspace{2mm}

In case the sum of the codimensions of two cycles is bigger than the
dimension of the space plus one, the definition of their algebraic distance 
is not as straightforward. On projective space however there are even several 
possibilities to define the algebraic definition in a useful way.

\satz{Definition}
For $p+q \geq t+1$ let 
$\Pe(W) \subset \Pe^t$ a subspace of codimension $q$, and $X$ an
effective cycle of pure codimension $p$ not meeting $\Pe(W)$ in $\Pe^t$.

\begin{enumerate}

\item
Let $G_W$ be the sub Grassmannian of the Grassmannian $\CG_{t+1,t+1-p}$
consisting of the subspaces of codimension $t+1-p$ that contain
$W$, and $V_XX$ as in Corollary \ref{VX}. Define
\[ D(\Pe(W),X) = D_G(\Pe(W),X) := \frac1{\deg G_W} D(G_W,V_X). \]
%\int_{G^{q'}_W} g_{V_X} \mu_G^{d_W} =
%   \int_{V \in G^{q'}_W} D(\Pe(V),X) \mu_G^{d_W} = \]
%\[  -\frac12 \int_{V \in G_W^{q'}} \int_{\Pe^t} (X|\Pe(V)) \mu_G^{d_W} = \]
%\[ -\frac12 \int_{V \in G_W^{q'}} \left( \int_X g_{\Pe(V)} -
%   \int_{\Pe^t} g_{\Pe(V)} \deg X \mu^{t+1-q'} \right) \mu_G^{d_W}, \]
%where $g_{\Pe(V)}$ is an admissible Green current for $\Pe(V)$, e.\@ g.\@
%$\Lambda_{\Pe(V)}$.

\item
In the same situation as in 1, define 
\[ D_\infty(\Pe(W),X) := {sup \atop V \in G_W} D(V,V_X) =
   {sup \atop V \in G_W} D(\Pe(V),X). \]
The equality holds by the Lemma \ref{dimension} and \ref{Dfunk1}.
%   -\frac12 {sup \atop V \in G_W^{q'}}  (X|\Pe(V)) = \]
%\[ -\frac12 {sup \atop V \in G_W^{q'}} \int_{\Pe^t} 
%   \left( \int_X g_{\Pe(V)} \mu^{t+1-q'} -
%   \int_{\Pe^t} g_{\Pe(V)} \deg X \mu^{t-p} \right). \]

\item
Let $Ch(X) \subset (\check{\Pe}^t)^{t+1-p}$ be the Chow divisor of $X$,
and 
\[ \left( {(t+1-p)(q-1) \atop q-1, \ldots, q-1} \right) =
   \frac{((t+1-p)(q-1))!}{(q-1)!^{t+1-q}} \]
the multinomial coefficients.
With $\Pe(\check{W})\subset \check{\Pe}^t$ the subspace dual to $\Pe(W)$,
define
\[ D_{Ch} (\Pe(W),X) := 
   \frac1{\left( {(t+1-p)(q-1) \atop q-1, \ldots, q-1} \right)}
   D(Ch(X),\Pe(\check{W})^{t+1-p}). \]
\end{enumerate}

\end{Satz}

%\satz{Remark} 
%If $p+q = t+1$, by Propositions \ref{glschnitt}, and \ref{Dfunk1}, 
%$D(\Pe(W),X) = D_\infty(\Pe(W),X) = D_0(\Pe(W),X)$.
%We will see later (Proposition \ref{abstrela}), that in this case also
%$D_{Ch}(\Pe(W),X) = D(\Pe(W),X) + c_7 \deg X$. 

%Note that the algebraic distances $D_0,D_\infty$ of two projective
%spaces are not necessarily symmetric, but will turn out to be symmetric
%modulo a constant.
%\end{Satz}

\satz{Fact} \label{addit}
For $p+q \geq t+1$, and $\Pe(W)$ a fixed subspace
of codimension $q$, the maps
\[ D_0(\Pe(W),\cdot),  D_{Ch}(\Pe(W),\cdot)  : Z_{eff}^p(\Pe^t) \to \R \]
are additive when defined.
\end{Satz}

\satz{Theorem} \label{egal}
\begin{enumerate}

\item
For $p+q\geq t+1$, any 
$X \in Z^p_{eff}(\Pe^t)$, and $\Pe(W)$ a subspace
of codimension $q$ not meeting $X$, 
\[ D_G(\Pe(W),X) \leq D_\infty(\Pe(W),X) \leq 
   c_3(t+1-q,t+1) \deg X + D_G(\Pe(W),X), \]

%\[ D(\Pe(F),X) = \log ||f_{V_X}|_{\Pe(F)}||_0 \leq
%   \log ||f_{V_X}|_{\Pe(F)}||_\infty \leq c' \deg X + D_0(\Pe(W),X). \]

\item
With $c_5,c_7$ the constants from Proposition \ref{chschieb} and
Lemma \ref{chnab}, for all $X$, $\Pe(W)$ as above,
\[ D_\infty(X,\Pe(W))-(c_7 + c_5) \deg X\leq D_{Ch}(X,\Pe(W)) \leq \]
\[ D_\infty(X,\Pe(W)) - c_5 \deg X. \]
\end{enumerate}
\end{Satz}

\proof
1. Let $f_X \in \Gamma(G,\CL^{\otimes \deg X})$.
be a  of norm $\log ||f_X||_0 = 0$ such that $V_X = \di(f_X)$. Then, 
\[ D_G(\Pe(W),X) = \frac1{\deg G_W} \int_{G_W} \log |f_X| \mu_G^{(p+q-t-1)(2t+2-p-q)}= 
\]
\[ \int_{G_W} \log |f_X| \omega_{G_W} = 
   \log ||f_X|_{G_W}||_0. \]
By Proposition \ref{normrel}, this is less or equal
\[ \log ||f_X|_{G_W}||_\infty = D_\infty(\Pe(W),X), \]
and this in turn by the same Proposition is less or equal
\[ \log ||f_X|_{G_W}||_0 + c_3(t+1-q,t+1) \deg X =
    D_G(\Pe(W),X) + c_3(t+1-q,t+1) \deg X. \]

\vspace{2mm}

2.
Let $f \in \Gamma((\check{\Pe}^t)^{t+1-p},O(\deg X,\ldots,\deg X))$ such
that $Ch(X) = \di f$. Then,
\[ |f_{[\check{v}_1,\ldots,\check{v}_{t+1-p}]}| = 
   \frac{f(\check{v}_1^{\deg X},\ldots,\check{v}_{t+1-p}^{\deg X})}
       {|\check{v}_1|^{\deg X} \cdots |\check{v}_{t+1-p}|^{\deg X}}, \]
and $f(\check{v}_1,\ldots,\check{v}_{t+1-p}) = 0$ if 
$\check{v}_i=\check{v}_j$ for some $i,j$. Hence
$|f_{[\check{v}_1,\ldots,\check{v}_{t+1-p}]}|$ takes it's supremum at some
$[\check{v}_1,\ldots,\check{v}_{t+1-p}]$ with $\check{v}_i \bot \check{v}_j$
for every $i \neq j$. 

Further, by Proposition \ref{chschieb} and Lemma \ref{Dfunk1},
\[ D(Ch(X),[\check{v}_1,\ldots,\check{v}_{t+1-p}]) =
   D(X^{t+1-p},\mbox{ker} \check{v}_1 \times 
   \cdots \times \mbox{ker} \check{v}_{t+1-p}), \]
which for $\check{v}_1,\ldots,\check{v}_{t+1-p}$ pairwise orthogonal
by Proposition \ref{chschieb} equals 
\[ D(X,\mbox{ker}{\check{v}_1} \cap \cdots \cap \mbox{ker}(\check{v}_{t+1-p}))-
   c_5(t,p,t+1-p)\deg X. \]
Hence,
\begin{equation} \label{c5}
\log |f|_{\Pe(\check{W})^{t+1-p}}|_\infty - \log |f_X|_0 = 
   D_{\infty}(\Pe(W),X)-c_5 \deg X. 
\end{equation}
Since
\[ D_{Ch}(\Pe(W),X) = \frac1{{(t+1-p)(q-1) \choose q-1,\ldots,q-1}}
   D(\Pe(\check{W})^{t+1-p},Ch(X)) = \]
\[ \frac1{{(t+1-p)(q-1) \choose q-1,\ldots,q-1}}
   \int_{\Pe(\check{W})^{t+1-p}} \log |f_X| \bar{\mu}^{(t-q)(t+1-p)} - \]
\[ \frac1{{(t+1-p)(q-1) \choose q-1,\ldots,q-1}} \int_{(\check{\Pe}^t)^{t+1-p}} 
   \log |f_X| \check{\mu}_1^q \cdots \check{\mu}_{t+1-p}^q\bar{\mu}^{(t-q)(t+1-p)}, \]
which by Lemma \ref{chnab} is less or equal
\[ \sup_{[\check{v}_1, \ldots,\check{v}_{t+1-p}] \in \check{W}^{t+1-p}}
   \log |(f_X)_{[\check{v}_1, \ldots,\check{v}_{t+1-p}]}| - \log |f_X|_0=
    \log |f_X|_{\Pe(\check{W})^{t+1-p}}|_\infty - \log |f_X|_0, \]
which in turn by the same proposition is less or equal
\[ \frac1{{(t+1-p)(q-1) \choose q-1,\ldots,q-1}}
   \int_{\Pe(\check{W})^{t+1-p}} \log |f_X| \bar{\mu}^{(t-q)(t+1-p)} - \]
\[ \frac1{{(t+1-p)(q-1) \choose q-1,\ldots,q-1}} \int_{(\check{\Pe}^t)^{t+1-p}} 
   \log |f_X| \check{\mu}_1^q \cdots \check{\mu}_{t+1-p}^q\bar{\mu}^{(t-q)(t+1-p)}+ 
   c_7 \deg X. \]
Together with equation (\ref{c5}), this implies the claim.

\subsection{Linear varieties}

For linear cycles, the algebraic distance is easy to calculate.
Let $\Pe(V), \Pe(W) \subset \C^{t+1}$ be properly intersecting
subspaces of dimension $p,q$, define 
$r:= \dim V \cap W -1= \max (-1,p+q-t)$, and let
$v_0, \ldots, v_{p+q+1-r} \in \C^{t+1}$ such that 
$(v_0,\ldots,v_r), (v_0,\ldots,v_p), \\
(v_0,\ldots,v_r,v_{p+1},\ldots,v_{p+q+1-r})$ are bases for 
$V\cap W, V, W$ respectively. Define
\[ |V,W| := \frac{\mbox{vol}\la v_0,\ldots,v_{p+q+1-r} \ra \;
                  \mbox{vol} \la v_0,\ldots,v_r \ra}
                 {\mbox{vol} \la v_0,\ldots,v_p \ra \;
                  \mbox{vol} \la v_0,\ldots,v_r,v_{p+1},\ldots,v_{p+q+1-r}\ra}. \]
Clearly, $\log |V,W| \leq 0$, and for $V,W$ both one dimensional
$|V,W|$ equals the Fubini-Study distance of the points $V,W$ in $\Pe^t(\C)$.

\satz{Proposition} \label{abst}
Let $\Pe(V),\Pe(W) \subset \Pe^t$ be projective subspaces of dimensions
$p,q$ respectively that intersect properly, and
with $r\leq t+1$ let $G_r$ be the Grassmannian  of $r$-dimensional
subspaces of $\C^{t+1}$. With the constants
constants $c_1,c_2$ and some positive constants $c_4,c_4$ only depending
on $p,q,r,t$, 
\begin{enumerate}
\item
For $p+q \geq t-1$,
\[ D(\Pe(V),\Pe(W)) = \log |V,W| + c_1(t,p,q) \leq
   c_1(t,p,q). \]
The inequality is an equality, iff the orthogonal complements of $V \cap W$ in
$V$ and $W$ are orthogonal to each other.

\item
There are positive constant $c_4,c_6$ only depending on $p,q,t$ such that for
$p+q<t$, 
\[ D_G(\Pe(V),\Pe(W))- c_4(p,q,t) = D_\infty(\Pe(V),\Pe(W))=\] 
\[ D_{Ch}(\Pe(V),\Pe(W)) - c_6(p,q,t) = \log |V,W| + c_2(t,q) \leq c_2(t,q). \]
The inequality is an equality iff $V,W$ are orthogonal to each other.
The algebraic distance of two projective spaces is thus symmetric modulo
$c_2(t,q) - c_2(t,p)+ c_4 (c_6)$.

\item
If $p+q+2 \leq r < t+1$, then $G_V,G_W$ in $G_r$ fulfill
\[ D(G_V,G_W) = \frac{\deg G_r}{\deg G_{V+W}}\log |V,W| + 
   c_8(t,p,q) \leq c_8(t,p,q). \]
The inequality is an equality, iff $V$ and $W$ are orthogonal to each other.

\item
For $p+q+1-t \geq r$, the varieties $G^F,G^E$ in $G_r$ fulfill
\[ D(G^F,G^E) = \frac{\deg G_r}{\deg G^{F \cap E}} D(\Pe(F), \Pe(E)) + c_8(t,p,q) 
   \leq c_8(t,p,q)+c_1(t,p,q). \]
The inequality is an equality iff the orthogonal complements of $E\cap F$
in $E$ and $F$ are orthogonal

\end{enumerate}
\end{Satz}

\proof
Let $v_0,\ldots,v_{p+q+1-r}$ as in the Definition of the distance of $V$ and $W$.
To shorten formulas, always assume that
the bases of $V\cap W, V$ and $W$ are orthonormal. Then,
\[ |V,W| = \mbox{vol} \la v_0,\ldots,v_{p+q+1-r} \ra. \]
Further, since distances and algebraic distances are invariant under the
unitary group, one may assume $v_i=e_i$ the standard unit vector 
for $i=0,\ldots,p$, hence
$\Pe(V)$ is the space where, $x_{p+1},\ldots,x_t$ vanish.

\vspace{2mm}

1. Use complete induction over the codimension of $\Pe(V)$. 
For $t-p = \mbox{codim} \; \Pe(V) = 1$, we have
$\log |V,W| = \log |v_t^\bot|=\log |(x_t)_{[v_t]}|$, where is $v_t^\bot$ 
projection of $v_t$ to
the orthogonal complement of $V=\la e_0,\ldots,e_{t-1} \ra$. By Lemma
\ref{normrelpr},
\begin{eqnarray*}
 D(\Pe(V),\Pe(W)) &=& \int_{\Pe(W)} \log |x_t| - \int_{\Pe^t} \log |x_t| \\
  &=& \sup_{P \in \Pe(W)} \log |(x_t)_P| - \sum_{j=1}^q \frac1j - 
   \int_{\Pe^t} \log |x_t| \\
   &=& \log |(x_t)_{[v_t]}|- \frac12 \sum_{j=1}^q \frac1j + 
   \frac12 \sum_{j=1}^t \frac1j,
\end{eqnarray*}
which by the above equals
\[ \log |V,W| + \sum_{j=q+1}^t \frac 1j =
   \log |V,W| + \sigma_t - \sigma_{t-1} - \sigma_q + \sigma_{q-1} =
   \log |V,W| + c_1(t,t-1,q). \]

Assume now the claim has been proved for 
$\mbox{codim} \; \Pe(V)= t-p$, and let \\
$\mbox{codim} \; \Pe(V)= t+1-p$, that is $\dim V =p$.
Again, because of the invariance of distance and algebraic distance under
the unitary group, one may assume that
$v_{p+1},\ldots,v_{t-1}$ are contained in $\la e_0,\ldots, e_{t-1}\ra$. Let
$\Pe(F)$ be the vanishing set of $x_t$, and $\Lambda_W$ as in (\ref{Levine}).
Then, by (\ref{Levint}),
\[ D(\Pe(V),\Pe(W)) = -\frac12\int_{\Pe(V)} \Lambda_{\Pe(W)} \mu^{t+1-p-q} - 
   +\sigma_{t-q} + \sigma_{q-1} - \sigma_t = \]
\[ -\frac12\int \Lambda_{\Pe(W)} \mu^{t+1-p-q} + 
   \frac12 \int_{\Pe(F)} \Lambda_{\Pe(W)} \mu^{t-p-q} -
   \frac12 \int_{\Pe(F)} \Lambda_{\Pe(W)} \mu^{t-p-q} +
   \sigma_{t-q} + \sigma_{q-1} - \sigma_t = \]
\[ -\frac12\int \Lambda_{\Pe(W)} \mu^{t+1-p-q} + 
   \frac12 \int_{\Pe(F)} \Lambda_{\Pe(W)} \mu^{t-p-q} + D(\Pe(F),\Pe(W)). \]
The first two summands are the algebraic distance of $\Pe(V)$ and $\Pe(W)$ 
as subvarieties of $\Pe(F)$. As $\Pe(V)$ has codimension $t-p$ in 
$\Pe(F)$, the induction hypothesis implies that this sum equals
\[ \log |\Pe(V),\Pe(W)\cap \Pe(F)|
   +  \sigma_{p+q-t-1}+ \sigma_{t-1} - \sigma_{p-1} -\sigma_{q-1} = \]
\[ \log \mbox{vol} \la v_0,\ldots,v_{t-1}\ra + 
   \sigma_{p+q-t-1}+ \sigma_{t-1} - \sigma_{p-1} -\sigma_{q-1}. \]
Further, as seen above
\[ D(\Pe(F),\Pe(W)) = \log |\Pe(F),\Pe(W)| + 
   \sigma_t-\sigma_{t-1}-\sigma_q+\sigma_{q-1} = \]
\[ \log |v_t^\bot| + \sigma_t-\sigma_{t-1}-\sigma_q+\sigma_{q-1}, \]
where $v_t^\bot$ denotes the projection of $v_t$ to the orthogonal complement
of $\Pe(F)$. Since 
$\mbox{vol}\la v_0,\ldots,v_{p+q+1-r}\ra =
\mbox{vol}\la v_0,\ldots,v_t\ra = \mbox{vol} \la v_0,\ldots,v_{t-1} \ra \cdot
|v_t^\bot|$, the claim follows with
$c_1(t,p,q) = \sigma_{p+q-t} + \sigma_t - \sigma_p - \sigma_q$.

\vspace{2mm}

2. By Definition and part one,
\[ D_{\infty} (\Pe(V),\Pe(W)) = \sup_{\Pe(F)} D(\Pe(F),\Pe(W)) =
   \sup_F \log |F,W| + c_1(t,t+1-q,q), \]
where $\Pe(F)$ runs over the superspaces of $\Pe(V)$ of codimension $q+1$.
As \\ $\sup_F \log |F,W| = \log |V,W|$, and the supremum is attained
at $F$ a direct sum of $W$ and an orthogonal complement
of $V+W$, the claim about $D_\infty$ follows with
$c_2(t,q) = c_1(t,t+1-q,q) = \sigma_1+\sigma_t-\sigma_{t+1-q}-\sigma_q$.

Next, let $f \in \Gamma(G_{t+1,t-q-1},L_G)$ be such that $V_{\Pe(W)} = \di f$.
This implies $f|_{G_V} \in \Gamma(G_V,L_{G_V})$. Since the fix group of $G_V$ 
inside the unitary group operates transitively on the unit sphere in
$\Gamma(G_V,L_{G_V})$, the number
\[ c_4 :=  \int_{G_V} \log |f| \omega_{G_V} - \sup_{P \in G_V} \log |f_P| \]
does not depend on $f,V,W$ but only on $t,p,q$, and we have
\[ D_G(\Pe(V),\Pe(W)) = \frac1{\deg G_V} D(G_V,V_{\Pe(W)}) = 
   \int_{G_V} \log |f| \omega_{G_V} - \int_G \log |f| \omega_G = \]
\[ c_4 + \sup_{P \in G_V} \log |f_P| - \int_G \log |f| \omega_G =
   c_4 + D_\infty(\Pe(V),\Pe(W)), \]
finishing the proof for $D_G(\bullet,\bullet)$. The claim about
$D_{Ch}(\bullet,\bullet)$ is proved similarly.

\vspace{2mm}

3. Assume first that $p+q+2=r$, and let $E$ be a subspace of
dimension $t+1-r$ that is orthogonal to $V+W$. Since $U(t+1)$ acts 
transitively on spaces $F,E$ of given dimensions with $F \bot E$,
\[ D(G_V,V_{\Pe(E)}) = c, \quad 
   D(G_V.G_W,V_{\Pe(E)}) = D(G_{V+W},V_{\Pe(E)})=c' \]
with constants $c,c'$ only depending on $p,q,r,t$. In 
Proposition \ref{doppelcor}.1 replace $p,q,p-q$ by $r,q+1,p+1$. Then,
by Lemma \ref{Dfunk1},
\[ D(G_V.V_{\Pe(E)},G_W) = 
   D(G_V(\varphi_r)_* \varphi_1^* \Pe(E),
      (\bar{\varphi}_r)_* \bar{\varphi}_{q+1}^* W) =
   D(\psi_{q+1})_* \psi_1^* \Pe(E)_V, W) = \]
\[ D(\Pe(E)_V, (\psi_1)_* \psi_{q+1}^* W) = D(\Pe(V+E), \Pe(W)) = \]
\[ \log |V+E,W| + c_1(t,t+p-r+1,q) = \log |V,W| +  c_1(t,t+p-r+1,q). \] 
The last inequality holding since $E$ is orthogonal to $V+W$. 
Hence, by Proposition \ref{vergl} and Lemma \ref{Dfunk1}
\begin{eqnarray*}
 D(G_V,G_W) &=& - D(G_V . G_W, V_{\Pe(E)}) + D(G_V . V_{\Pe(E)},G_W) + 
   \frac1{\deg G_W} D(G_V,V_{\Pe(E)}) \\
 &=&  - c' + \log |V,W| +  c_1(t,t+p-r+1,q) + \frac{c}{\deg G_W}, 
\end{eqnarray*}
and the claim holds with 
$c_8(t,p,q) = c'_8(t,p,q):=- c' + c_1(t,t+p-r+1,q)+\frac{c}{\deg G_W}$.

Let now $p+q+2 < r$ and $E$ a space of dimension $r-p-q-1$ that is
orthogonal to $V+W$. Then, using again a transitive action of $U(t+1)$, one
gets $D(G_V,G_E) =c''$ with $c''$ only depending on 
$p,t$. Further, by the first half of the proof
\[ D(G_{V+W},G_E) = \log |V+W,E| + c_8'(t,p+q+1,r-p-q-2) =
   c_8'(t,p+q+1,r-p-q-2), \]
and
\[ D(G_{V+E}, G_W) = \log |V+E,W| + c_8'(t,t-q,q) = \log |V,W| + c_8'(t,t-q,q). \]
Hence, by Proposition \ref{vergl},
\[ \frac{\deg G_{V+W}}{\deg G_r} D(G_V,G_W) = \frac1{\deg G_E} D(G_V,G_W) = \]
\[ - D(G_V.G_W,G_E) + \frac1{\deg G_W} D(G_V,G_E) + D(G_V.G_E,G_W) = \]
\[  -D(G_{V+W}, G_E) + D(G_V,G_E) + D(G_{V+E}, G_W)  = \]
\[ c_8'(t,p+q+1,t-p-q-2) + \frac{c''}{\deg G_W} +c_8'(t,t-q,q) + 
   \log |V,W|, \]
finishing the proof.

\vspace{2mm}

4. Follows similarly to part 3, this time using the orthogonal complement
$E$ of $V \cap W$ and Proposition \ref{doppelcor}.2.

\subsection{Decompositions}

With $W$ a projective K\"ahler variety of dimension $t$ and
$Z \subset W$ a subvariety of codimension $r$, and  
$X,Y$ effective cycles of pure codimensions $p,q$, we have the algebraic 
distance $D(X,Y)$ if $p+q\leq t+1$ and $X,Y$ intersect properly.
If for $p+q \leq t+r+1$ the supports of $X$ and $Y$ are contained in $Z$, and
they intersect properly as cycles in $Z$, denote with $D^Z(X,Y)$ the algebraic 
distance of $X,Y$ as cycles in $Z$.

\satz{Proposition} \label{raumraum}
Let $t,p,r\leq q$ be natural numbers with $p+r \leq t+1$, further
$X \in Z_{eff}^p(\Pe^t)$ and $\Pe(W) \subset \Pe(F) \subset \Pe^t$ subspaces
of codimensions $q,r$ such that the intersection of $\Pe(W)$ and
$\Pe(F)$ with $X$ are proper. 

\begin{enumerate}

\item
If $p+q \leq t+1$, then
\[ D(\Pe(W),X) = D^{\Pe(F)}(X.\Pe(F),\Pe(W)) + D(\Pe(F),X). \]

\item
If $p+q >t+1$, and $X_W$ is the cycle defined in (\ref{XW}), then
\[ D(\Pe(W),X) = \]
\[ D^{\Pe(F)}(\Pe(W),X.\Pe(F)) + D(X,\Pe(F))-D(X_W,\Pe(F)) +
   c_{11} \deg X, \]
where $c_{11}$ is a constant only depending on $p,q$, and $t$.

\item
If $p+q \leq t+1$, then
in the Grassmannian $G = G_{p,t+1}$,
\[ D(V_X,G^W) = \frac{\deg G^W}{\deg G^F} \;D(V_X, G^F) + D^{G^F}(V_X.G^F,G^W). \]

\item 
If $p+q \leq t+1$, then in the Grassmannian  $G=G_{p,t+1}$,
\[ D(V_X,G^F) = \deg G^F D(X,\Pe(F)) \leq d(t-q,p) c_2(t,p,q) \deg X, \]
with $d$ a constant depending on $t-q$ and $p$.

\end{enumerate}
\end{Satz}

\proof
1. Let $g_X$ be an admissible Green form for $X$. By \cite{SABK}, Lemma II.2,
$g_X^{\Pe(F)} := g_X|_{\Pe(F)}$ is a Green form of log type for 
$X.\Pe(F)$ in $\Pe(F)$, which by
Proposition \ref{funk} is admissible. Hence,
\[ D(X,\Pe(W)) = -\frac12 \int_{\Pe(W)} g_X \mu^{t+1-p-q} + 
   \frac12 \int_{\Pe^t} g_X \mu^{t+1-p} = \]
\[ -\frac12 \int_{\Pe(W)} g_X^{\Pe(F)} \mu^{t+1-p-q} + \frac12 \int_{\Pe(F)} 
   g_X^{\Pe(F)} \mu^{t+1-p-r} -
   \frac12 \int_{\Pe(F)} g_X \mu^{t+1-p-r} + \]
\[ \frac12 \int_{\Pe^t} g_X \mu^{t+1-p} = 
   D^{\Pe(F)}(X.\Pe(F),\Pe(W)) + D(X,\Pe(F)). \]

\vspace{2mm}

2. If $\dim X =0$, then $\Pe(F) = \Pe^t$, hence the equality is trivial.

For $\dim X > 0$, let $V \subset F$ be a subspace of dimension 
$p+q-t-1$ such that 
$U:= V+W$ is a direct sum and $\Pe(U)$ does not meet the support of $X$.
One has $U = G_V \cap G_W$ in $G = G_{t+1,p}$. With $V_X$ as in 
Corollary \ref{VX}, one has $H(\delta_{V_X}) = \deg X \mu_G$, and by Proposition
\ref{vergl}, with $g_W$ a normalized Green form for $G_W$, and 
$\eta_W=H(\delta_{G_W})$,
\[ \deg X D(G_V,G_W)+ D(V_X,U) = \]
\begin{equation} \label{Punktraum00}
-\frac12 \deg X \int_{G_V} g_W \mu_G + D(V_X, U) =
   D(G_W. V_X,G_V) - \frac12 \int_{V_X} g_W \; \eta_{G_V}. 
\end{equation}
If $f \in \Gamma(G,L^{\otimes \deg X})$ with $V_X = \di f$, the associativity
of the star product implies
\[ -\frac12 \int_{V_X} g_W \eta_{G_V} =
   - \frac12 \deg X \int_G g_W \eta_{G_V} \mu_G +
   \int_{G_W} \log |f| \; \eta_{G_V} - \int_G \log |f| \; \eta_{G_W} \eta_{G_V} \]
By Proposition \ref{gras}, $\eta_{G_V} \eta_{G_W} = \omega_G$, and as
$g_W$ is normalized, this equals
\[ \frac1{\deg G_W} D(G_W,V_X). \]

Inserting into (\ref{Punktraum00}) gives
\[ \deg X D(G_V,G_W) + D(V_X,U) = D(G_W.V_X,G_V) + 
   \frac1{\deg G_W} D(G_W,V_X). \]
Thus, by (\ref{GWgl}), Proposition \ref{doppelcor}, 
Proposition \ref{Dfunk1}, and the Definition of $D(\Pe(W),X)$,
\begin{equation} \label{Punktraum01}
\deg X D(G_W,G_V) + D(X,\Pe(U)) =
   D(X_W,\Pe(V)) + D(\Pe(W),X). 
\end{equation}
By part one,
\[ D(X,\Pe(U)) = D^{\Pe(F)}(X.\Pe(F), \Pe(U)) + D(X,\Pe(F)), \]
and 
\[ D(X_W,\Pe(V)) = D^{\Pe(F)}(X_W.\Pe(F),\Pe(V)) + D(X_W,\Pe(F)),\]
hence
\[ \deg X D(G_W,G_V) + D^{\Pe(F)}(X.\Pe(F),\Pe(U)) + D(X,\Pe(F)) = \]
\begin{equation} \label{Punktraum02}
D^{\Pe(F)}(X_W.\Pe(F),\Pe(V)) + D(X_W,\Pe(F)) + D(\Pe(W),X). 
\end{equation}

Further repeating the argument for (\ref{Punktraum01}) 
with $\Pe(F)$ instead of $\Pe^t$ and
$X. \Pe(F)$ instead of $X$, one gets
\[ \deg X D^{G^F}(G_V^F,G_W^F) + D^{\Pe(F)}(X.\Pe(F),\Pe(U)) = \]
\begin{equation} \label{Punktraum11}
D^{\Pe(F)}((X.\Pe(F))_W,\Pe(V)) + D^{\Pe(F)}(\Pe(W), X.\Pe(F)).
\end{equation}

Using $(X.\Pe(F))_W = X_W . \Pe(F)$ while
subtracting (\ref{Punktraum02}) and (\ref{Punktraum11}) leads to
\[ \deg X \left( D(G_V,G_W) - D^{G^F}(G_V^F,G_W^F)\right) + D(X,\Pe(F)) = \]
\[ D(X_W,\Pe(F))  + D(\Pe(W),X) - D^{\Pe(F)}(\Pe(W),X.\Pe(F)). \]

Finally, by Proposition \ref{abst}.3,
\[ D(G_V,G_W) - c_8(t,p,q) = D^{G^F}(G_V^F,G_W^F) -c_8(p,p,q) =
   \log |\Pe(V),\Pe(W)|, \]
hence the claim follows with $c_{11}(t,p,q)=c_8(p,p,q)-c_8(t,p,q)$.

\vspace{2mm}

3. Let $f \in \Gamma(G,L^{\otimes \deg X})$ such that $V_X= \di f$, and
$\eta_{G^W} = H(\delta_{G^W})$. Then,
\[ D(V_X,G^W) = \int_{G^W} \log |f| \mu_G^{p(t+1-q-p)} -
               \int_G \log |f| \eta_{G^W} \mu_G^{p(t+1-q-p)}, \]
which by Proposition \ref{gras}.5 equals
\[ \int_{G^W} \log |f| \mu_G^{p(t+1-q-p)} -
   \deg G^W \int_G \log |f| \omega_G = \]
\[ \int_{G^W} \log |f| \mu_G^{p(t+1-r-p)} -
   \int_{G^F} \log |f| \eta_{G^W}^{G^F} \mu_G^{p(t+1-q-p)} + \]
\[ \int_{G^F} \log |f| \eta_{G^W}^{G^F} \mu_G^{p(t+1-q-p)}
   - \deg G^W \int_G \log |f| \omega_G, \]
where $\eta_{G^W}^{G^F} = H^{G^F}(\delta_{G^W})$ with $H^{G^F}$ the harmonic 
projection in the sub Grassmannian $G^F \subset G$.
Using Proposition \ref{gras}.5 once more, this equals
\[ \int_{G^W} \log |f| \mu_G^{p(t+1-q-p)} -
   \int_{G^F} \log |f| \eta_{G^W}^{G^F} \mu_G^{p(t+1-q-p)} + \]
\[ \frac{\deg G^W}{\deg G^F} \int_{G^F} \log |f| \mu_G^{p(t+1-r-p)} -
   \frac{\deg G^W}{\deg G^F} \int_G \log |f| \eta_{G^F} \mu_G^{p(t+1-r-p)} = \]
\[ D^{G^F}(V_X .G^F, G^W) + \frac{\deg G^W}{\deg G^F} D(V_X,G^F). \]

\vspace{2mm}

4. Assume first, $\Pe(F) = \Pe(W)$ has codimension $t+1-p$. Then,
$G^W=G_W=W$ is a point, hence has degree one, and
\begin{equation} \label{compldimschieb}
D(X,\Pe(W)) = D(V_X,W) = D(V_X,G_W) = D(V_X,G^W)
\end{equation} 
follows from Lemma \ref{dimension} and Lemma \ref{Dfunk1}, and 
Remark \ref{remark}.

In the general case, let $\Pe(V) \subset \Pe(F)$ be a subspace of
dimension $p-1$ that does not meet the support of $X$. 
By part one,
\[ D(X,\Pe(F)) = D(X,\Pe(V)) - D^{\Pe(F)}(X.\Pe(F),\Pe(V)). \]
Thus, applying (\ref{compldimschieb}) once for $\Pe^t$ and $G_{t+1,}$ and once
for $\Pe(F)$ and $G^F$,
\[ D(X,\Pe(F)) = D(V_X,V) - D^{G_F}(V_{X.\Pe(F)},F) =
   D(V_X,V) - D^{G_F}(V_X.G^F,W), \]
which by part 3 equals 
\[ D(V_X,W) -\left(-\frac1{\deg G^F} D(G^F,V_X)+D(V_X,W)\right) = 
   \frac1{\deg G^F} D(G^F,V_X). \]
Since, by Remark \ref{remark}, $D(X,\Pe(F)) \leq c_2(t,p,q)$, the
claim follows with $d(t-q,p) = \deg G^F$.

\subsection{Reduction to distances to points}

\satz{Lemma} \label{linausw}
For $p,q \leq t+1$ let $X \in Z_{eff}^p(\Pe^t)$, and
$\Pe(V) \supset \Pe^t$ a subspace of codimension $q$
intersecting $X$ properly, and $s$ some natural number between $0$ and $t+1$.
There is a positive constant $C$, only depending on $t,p,q,s$ such that,

\begin{enumerate}
\item
if $p+q=t+1$, then
\[ D(\Pe(E),X) \geq D(\Pe(V),X) - C(t,p,q,s) \deg X \]
for every space $\Pe(E)$ of dimension $s$ such that 
$\Pe(E) \subset \Pe(V)$ or $\Pe(E) \supset \Pe(V)$.

\item
if $p+q < t+1$, then
\[ D(\Pe(E),X) \geq D(\Pe(V),X) - C(t,p,q,s) \deg X, \]
for every subspace $\Pe(E)$ of dimension $s>t-q$ containing $\Pe(V)$, and
\[ \sup_{E \in G_s^V} D(\Pe(E),X) \geq D(\Pe(V),X)- C(t,p,q,s) \deg X, \]
for every $s \leq t-q$.

\item
if $p+q > t+1$, then
\[ D(\Pe(E),X) \geq D(\Pe(V),X)- C(t,p,q,s) \deg X, \]
for every subspace $\Pe(E) \subset \Pe(V)$ of dimension $s<t-q$, and
\[ \sup_{E \in (G_s)_V} D(\Pe(E),X) \geq D(\Pe(V),X) - C(t,p,q,s) \deg X, \]
for every $s \geq t-q$.
\end{enumerate}
\end{Satz}

\proof
Let $f \in \Gamma(G_p,L_G^{\otimes \deg X})$ such that $V_X = \di f$.

\vspace{2mm}

1. Let first $\Pe(E)$ be a subspace of dimension $s \leq t+1-q$ such that
$\Pe(E) \subset  \Pe(V)$. Then, with $d(s,p,t) = \deg G_E$,
\[ D(\Pe(E),X) = \frac1{\deg G_E} D(G_E,V_X) = \]
\[ \frac1{d(s,p,t)}\int_{G_E} \log |f| \mu^{(p-s)(t+1-p)} - 
   \frac1{d(s,p,t)}\int_{G_p} \log |f| \eta_{G_E} \mu^{(p-s)(t+1-p)}, \]
which by Propositions \ref{normrel} and \ref{gras}.5
is greater or equal
\[ \log |f_V| - \int_{G_p} \log |f| \; \omega_{G_p} -c_3(s,t+1) \deg X= 
   D(V,V_X) -c_3(s,t+1) \deg X, \]
which by Proposition \ref{Dfunk1} equals 
\[ D(\Pe(V),X)-c_3(s,t+1) \deg X. \]

If $s > t+1-q$ and $\Pe(E) \supset \Pe(V)$, then by Proposition 
\ref{raumraum}.1,
\[ D(\Pe(E),X) = D(\Pe(V),X)- D^{\Pe(E)}(\Pe(V),X. \Pe(E)), \]
which by Remark \ref{remark} is greater or equal
\[ D(\Pe(V),X) - \bar{c}_1(t-s,t+1-p-s,p) \deg X. \]

\vspace{2mm}

2.
%Let $p+q\leq t+1$ and $\Pe(F) \supset \Pe(V)$ a space of codimension $r\leq q$.
%By Proposition \ref{raumraum}.2,
%\[ D(\Pe(F),X) = - D^{\Pe(F)}(\Pe(V),X) + D(\Pe(V),X), \]
%which by Ramark \ref{remark} is greater or equal
%$D(\Pe(V),X) - c_1(t,p,q) \deg X$. 
%
By Proposition \ref{raumraum}.4, with $d(p,q,t) = \deg G^V$,
\[ D(X,\Pe(V)) =  \frac1{\deg G^V} \; D(V_X,G^V)= \]
\[ \frac1{d(p,q,t)} \int_{G_p^V} \log |f| \; \mu_G^{p(t+1-q-p)}- 
   \frac1{d(p,q,t)} \int_{G_p} \log |f| \; \mu_G^{p(t+1-p)}, \]
which by Propositions \ref{normrel} and \ref{gras}.5 is less or equal
\[ \sup_{W \in G_p^V} \log|f_W|- 
   \frac1{d(t,q)} \int_{G_p} \log |f| \mu_G^{p(t+1-p)} = 
   \sup_{W \in G_p^V} \log|f_W| -
   \int_{G_p} \log |f| \; \omega_{G_p} = \]
\[ \sup_{W \in G_p^V} D(V_X,W) =
   \sup_{W \in G_p^V} D(\Pe(W),X). \]
Let $W_0 \in G_p^V$ be such that $\sup_{W\in G_p^V} D(\Pe(W),X) = D(\Pe(W_0),X)$.
Then, \\ $D(\Pe(W_0),X) \geq D(\Pe(V),X)$, any space containing
$\Pe(V)$ also contains $\Pe(W_0)$, and for every $s \leq t-q$, there is a space
$\Pe(E)$ of codimension $s$ that either contains $\Pe(W_0)$ or is contained
in $\Pe(W_0)$. Thus, the claim follows from part one.

\vspace{2mm}

3. By definition
\[ D(\Pe(V),X) = \frac1{\deg G_V} D(V_X,G_V) = \]
\[ \frac1{\deg G_V} \int_{G_V} \log |f| \mu_G^{(p-q)(t+1-p)} - 
   \frac1{\deg G_V} \int_{G_p} \log |f| \eta_{G_V} \mu_G^{p(t+1-p)}, \]
which by Proposition \ref{normrel} is less or equal
\[ \sup_{W \in (G_p)_V} \log |f_W| -
   \frac1{\deg G_V} \int{G_p} \log |f|  \eta_{G_V} \mu_G^{(p-q)(t+1-p)} = \]
\[ \sup_{W \in (G_p)_V} D(V_X,W) = \sup_{W \in G_V} D(X,\Pe(W)). \]
Choose $W_0$ of dimension $p$ such that 
$D(X,\Pe(W_0)) = \sup_W D(X,\Pe(W))$.
Then, $D(\Pe(W_0).X) \geq D(\Pe(V),X)$, any space contained in
$\Pe(V)$ is also contained in $\Pe(W_0)$, and for every $s \geq q$ there is
a space $\Pe(E)$ of codimension $s$ that is either contained in $\Pe(W_0)$
or contains $\Pe(W_0)$. Thus the claim follows again from part 1.

\satz{Definition} \label{Dpt}
For $p+q \geq t+1$, 
$X \in Z^p_{eff}(\Pe^t_\C)$, and 
$\Pe(W) \subset \Pe^t_\C$ a subspace of codimension $q$ not meeting $X$, let
$X_W$ be the cycle defined in (\ref{XW}), and $\Pe(F) \subset \Pe^t$  a
subspace of dimension at least $p$ containing
$\Pe(W)$ and intersecting $X$ properly.

\begin{enumerate}
\item
The space $\Pe(F)$ is called $c$-admissible for $X$ and $\Pe(W)$ if
\[ D(\Pe(F),X) \geq -c \deg X \quad \mbox{and} \quad
   D(\Pe(F),X_W) \geq -c \deg X. \]
\item
If the dimension of $\Pe(F)$ equals $p$, define
\[ D_{\Pe(F)}(\Pe(W),X) := D^{\Pe(F)}(\Pe(W),X.\Pe(F)) - (c_2+c_4) \deg X = \]
\[ \sum_{x \in supp(\Pe(F).X)} n_x \log |\theta,x|, \]
where the $n_x$ are the intersection multiplicities of $\Pe(F)$ and $X$ at $x$.
%and
%$\Pe(V) \subset \Pe(\bar{V}) \subset \Pe^t$ subspaces of
%dimensions $p+q-t-1$ and $p+q-t$ respectively, such that the function
%\[ \CF_{(p+q-t-1,p+q-t)} \to \R, \quad (V', \bar{V}') \mapsto
%   D(X,\Pe(\bar{V}')) + D(X_W, \Pe(V')) \]
%takes its supremum at $(V,\bar{V})$. Let further $\Pe(F)\subset \Pe^t$ be the
%subspace of dimension $p+1$ spanned by $\Pe(W)$ and $\bar{V}$, and define
\end{enumerate}
\end{Satz}

%By fact \ref{abst}, 
%\begin{eqnarray*}
%c_1 \deg X + D_{pt}(\theta,X) &=& 
%     D_0^{\Pe(F)}(\theta,X . \Pe(F))- \deg X c_4(p,q,t) \\
%&=& D^{\Pe(F)}_\infty(\theta,X . \Pe(F))- \deg X c_2(t,q,t+1-q) \\
%&=& D^{\Pe(F)}_{Ch}(\theta,X . \Pe(F)) - \deg X c_6(t,q,t+1-q), 
%\end{eqnarray*}
%where $\Pe(F)$ is the subspace minimizing $D_{\Pe(F)}(\theta,X)$.

%\satz{Theorem} \label{Punktraum}
%Let $X \in Z^p_{eff}(\Pe^t_\C)$, and $\Pe(W) \in \Pe^t$ a subspace of
%codimension $q>t-p$ not meeting $X$. 
%\begin{enumerate}

%\item
%If $\Pe(F)$ is a subspace
%of codimension $t-p \leq r<q$ containing $\Pe(W)$ such that the intersection 
%$\Pe(F) . X$ is proper and equals a sum of projective spaces. Then,
%\[ D_{Ch}(\Pe(W),X) \leq \sum_{x \in supp(\Pe(F). X)}n_x \log |\Pe(W),x|,\]
%where $n_x$ denotes the intersection multiplicity. The condition that
%$\Pe(F) . X$ is a sum of projective spaces is of course fullfilled if
%$r = t-p$.

%\item
%If in 1, $r = t-p$, and
%$D(\Pe(F),X)$ is maximal among all $D(\Pe(F'),X)$ with
%$\Pe(F')$ a subspace of codimsion $t-p$ containing $\Pe(W)$, then
%\[ \sum_{x \in supp(\Pe(F). X)} n_x \log |\Pe(W),x| \leq
%   D_\infty(\Pe(W),X) + c \deg X. \]
%\end{enumerate}
%\end{Satz}

%Theorem \ref{egal}, and \ref{Punktraum} together imply

\satz{Theorem} \label{Punktraum}
With the above notations,
\begin{enumerate}
\item
let $X \in Z^p(\Pe^t)$ and $\Pe(W)$ a subspace of codimension
$q > t-p$ not meeting $X$, and $c_{10}(j,t+1)) \sum_{j=1,j\neq i}^t c_3(j,t+1)$. 
For every $s \geq t-p$, and every 
$c\geq c_{10}+C(p,q,s,t)$, with $C$ the constant from the previous Lemma,
there is at least one space $\Pe(F)$ of dimension $r$ that is $c$-admissible
for $X$ and $\theta$.

\item
There are positive constants $e_1,e_2,e_1',e_2'$ such that for all $p \leq t$
and every $X \in Z^p_{eff}(\Pe^t)$ and
$\theta \in  \Pe^t$ not contained in the support in $X$,
if $\Pe(F)$ is a $c$-admissible subspace of dimension 
$q \geq p$ for $X$ and $\Pe(W)$, then
\[ D(\theta,X) \leq D^{\Pe(F)} (\theta,X. \Pe(F)) + (c+e'_1) \deg X \leq 
   D(\theta,X) + (2c+e_1'+e_2') \deg X. \]
In particular, if $q=p$,
\[ D(\theta,X) \leq D_{\Pe(F)} (\theta,X) + c+e_1 \deg X \leq 
   D(\theta,X) + (2c+e_1+e_2) \deg X. \]
Hence, the algebraic distance of $\theta$ to $X$ essentially equals
the weighted sum of the distances of $\theta$ to the points contained 
in the intersection of $X$ with some projective subspace
$\Pe(F) \subset \Pe^t$ of codimension $t-p$ containing $\theta$.
\end{enumerate}
\end{Satz}

\vspace{2mm}

{\bf Remark}
If $X$ is a hypersurface, the stronger estimate
\[ D(\theta,X) \leq 
   \inf_{\Pe(F)\ni x,\dim \Pe(F)=1} D_{\Pe(F)}(\theta,X) + e_3 \deg X, \]
holds. This is wrong in general as the following example shows:
Let 
\[ t=3, \quad \theta = [1,0,0,0], \quad 
   X=\di(x^{D-1}w-z^D) . \di(\epsilon x-y), \quad \Pe(F) = \di(w), \]
\[ \Pe(V) = \di (y-w) . \di(z-w), \quad \Pe(E) = \di(z-w), \]
\[ Y=\di(z-w).\di(x^{D-1}w-z^D), \]
and $\zeta_{D-1}$ a primitive $(D-1)$th root of unity. 
By Lemmas  \ref{normrelpr} and \ref{linausw}.1, since
$[1,0,1,1]$ lies in the support of $\di z-w)$,
\[ D(\di(z-w),\di(x^{D-1}y-z^D)) \geq D([1,0,1,1],\di(x^{D-1}y-z^D))-CD = \]
\[ \log |1^{D-1} \cdot 0-1^D| - \log |x^{D-1}y-z^D|_0 -C D \geq \]
\[ -\log |x^{D-1}y-z^D|_\infty- D \left(C+ \frac12 \sum_{i=1}^t \frac1i\right) \geq
   - D \left(C+ \frac12 \sum_{i=1}^t \frac1i\right). \]
Further, with $\Pe(\bar{E}) = \di \; y$ and $\pi$  the orthogonal projection
$\pi:\Pe^3 \setminus [0,1,0,0] \to \Pe(\bar {E})$, one has
$Y = \overline{\pi^*(Y.\Pe(\bar{E}))}$, thus by \cite{BGS}, Proposition
5.1.1, $D(Y,\Pe(\bar{E})) = c_1D>0$, and consequently
$D(Y,\di(\epsilon x-y)) \geq 0$ for $\epsilon$ sufficiently small. 
Hence, by Proposition \ref{vergl} and Remark \ref{remark}
\[ D(\Pe(E),X) = -D(\di(x^{D-1}w-z^D),\di (\epsilon x-y)) +   \]
\[ D(\di(z-w),\di(x^{D-1}w-z^D))+ D(Y,\di(\epsilon x-y)) \geq  
   -D \left(C+\bar{c}_1+ \frac12 \sum_{i=1}^t \frac1i\right). \]
for $\epsilon$ small, and since $\theta \in \Pe(V)$, by Proposition
\ref{raumraum}.1, and \ref{abst}.2
\[ D_{\infty}(\theta,X) \geq D(\Pe(V),X) = 
   D(\Pe(E),X) + D^{\Pe(E)}(\Pe(V),\Pe(E).X) = \]
\[ D(\Pe(E),X) + 
   D^{\Pe(E)}\left([1,\epsilon,0,0]+\sum_{i=1}^{D-1} 
      [1,\epsilon,\zeta_{D-1}^i,\zeta_{D-1}^i], \Pe(V)\right) \geq \]
\[ D(\Pe(E),X) + \log \left| [1,\epsilon,0,0]+\sum_{i=1}^{D-1} 
      [1,\epsilon,\zeta_{D-1}^i,\zeta_{D-1}^i], \Pe(V)\right| +(c_2+c_4)D \geq \]
\[ -D\left(C+ \bar{c}_1 + c_2 + c_4+\frac12 \sum_{i=1}^t \frac1i\right) + 
   \log \frac{\epsilon}{\sqrt{1+\epsilon^2}} + (D-1) \log \frac12  \geq
   \log \epsilon -cD, \]
for small $\epsilon$, and some fixed constant $c$.

On the other hand, $X. \Pe(F) =D [1,\epsilon,0,0]$, hence
\[ D_{\Pe(F)}(\theta,X) = D \log |[1,0,0,0],[1,\epsilon,0,0]| =
   D \log \frac{\epsilon}{\sqrt{1+\epsilon^2}} \leq
   D \log \epsilon, \]
and the inequality 
\[ D_\infty(\theta,X) \leq 
   \inf_{\Pe(F)\ni x,\dim \Pe(F)=1} D_{\Pe(F)}(\theta,X) + c' \deg X \]
for some fixed $c'$ would imply
\[ \log \epsilon -c D \leq D \log \epsilon + c' D, \]
which is clearly wrong for $\epsilon$ sufficiently small and $D \geq 2$.

%\satz{Corollary} \label{Suk}
%Let $X \in Z^p(\Pe^t_\C)$, and $\Pe(W) \in \Pe^t$ a subspace of
%codimension $r>t-p$ not meeting $X$; further
%$\Pe(F)$ a superspace of dimension
%$p$ of $\Pe(V)$. Then,
%\[ D_\infty(\Pe(W),X) =
%   D(X,\Pe(F)) + \sum_{x \in X . \Pe(F)}
%   n_x \log |\Pe(W),x| + c(p,p,r,t) \deg X \]
%where $n_x$ are the intersection multiplicities.
%\end{Satz}

%\proof
%Follows from the Proposition and Corollary \ref{abst}.

\vspace{2mm}

\satz{Proposition}
For $i=1,\ldots,t$ let $X_i \in Z^1_{eff}(G_i)$ be an effective 
cycle of codimension one in the corresponding Grassmannian over
$\C$. There exists a flag of vector spaces
\[ \{0\} \subset V_1 \subset \cdots \subset V_t \subset \C^{t+1} \]
with $\dim V_i =i$ such that
\[ D(X_i,V_i) \geq -c_{10}(i) \deg X_i, \]
where $c_{10}(i) := \sum_{j \neq i} c_3(j)$.
\end{Satz}

\proof
Let $\bar{X}_i$ be the cycle 
$\bar{X}_i = \left(\sum_{k \neq i} \deg X_k\right) \cdot \; X_i$. Then,
$D:= \prod_{k=1}^t \deg X_k= \deg \bar{X}_i$.
Next, let $F$ be the complete flag variety, $\varphi_i:F \to G_i,i=1,\ldots,t$ 
the canonical projections, and $Y_i := \varphi_i^* \bar{X}_i$. Then,
\[ \left[\sum_{i=1}^t Y_i \right] = c_1(L^{\otimes D}), \]
with $L$ the line bundle on $F$ introduced in section 3.
Let $f_i \in \Gamma(G_i, L_i^{\otimes D})$ be such that $\bar{X}_i = \di f_i$.
Also assume that $f_i$ is normalized, i.\@ e.\@ 
$\int_{G_i} \log |f_i| \omega_i = 0$.
With $f = \prod_{i=1}^t \varphi_i^* f_i$, we have
$f \in \Gamma(F,L^{\otimes D})$ and $\sum_{i=1}^t Y_i = \di f$.
Further, $\sum_{i=1}^t \varphi_i^* (\log |f_i|) = \log |f|$, hence
\[ \int_F \log |f| \omega_F = 
   \sum_{i=1}^t \int_F \varphi_i^* (\log |f_i|) \omega_F =
   \sum_{i=1}^t \int_{G_i} \log |f_i| (\varphi_i)_* \omega_F =
    \sum_{i=1}^{t+1} 0. \]
The last equality holds because $(\varphi_i)_* \omega_F= \omega_i$ since
$(\varphi_i)_* \omega_F$ is $U(t+1)$-invariant and
$\int_{G_i} \omega_i = \int_{\varphi^*_i G_i} \omega_F=
\int_{G_i} (\varphi_i) \omega_i = \int_F \omega_F =1$.
By Proposition \ref{normrel}, there is a point $P \in F$ such that
\[ \log |f_P| \geq \int_F \log |f| \omega_F = 0, \]
and
\[ \log |(f_i)_{\varphi_i(P)}| \leq \log |f_i|_\infty \leq 
   c_3(i) \; D + \log |f_i|_0 = c_3(i) \; D, \quad i=1,\ldots,t+1,  \]
hence, 
\[ -\int_F \varphi_i^*(\log |f_i|) \omega_F + \log |(\varphi_i^* (f_i))_P| = \]
\[ \sum_{j \neq i} \int_F \varphi_j^*(\log |f_j|) \omega_F - 
   \int_F \log |f| \omega_F + 
   \log|f_P| - \sum_{j\neq i} \log |(\varphi_j^*(f_j))_P| = \]
\[ \sum_{j \neq i} 
   \left(\int_{G_j} \log |f_j| \omega_{G_j}-|(f_j)_{\varphi_j(P)}| \right) +\log |f_P|
   \geq - D \sum_{j \neq i} c_3(j) = -c_{10}(i) D. \] 
The point $P$ corresponds to a complete flag
$\{0\} \subset V_1 \subset \cdots \subset V_t \subset \C^{t+1}$ with
$V_i = \varphi_i(P)$. Thus,
\[ D(\bar{X}_i,V_i) = \log|(f_i)_{V_i}|-\int_{G_i} \log |f_i| \omega_i =
    \int_{G_i} \log |f_i| \delta_{(\varphi_i)_* P} - 
    \int_{G_i} \log |f_i|  (\varphi_i)_* \omega_F = \]
\[ \int_F \varphi_i^* (\log |f_i|) \delta_P - 
   \int_F \varphi_i^* (\log |f_i|) \omega_F \geq - c_{10}(i) D. \]
Dividing the inequality by $\prod_{k \neq i} \deg X_i$ gives
\[ D(X_i,V_i) \geq - c_{10}(i) \deg X_i \]
as claimed.

\satz{Corollary} \label{adm}
Let $X \in Z_{eff}(\Pe^t)$ be an effective cycle in $\Pe^t$, and
denote by $X_i$ its component of codimension $i$, $0=1,\ldots,t$. There
is a flag of vector spaces
\[ \{0\} \subset V_1 \subset \cdots \subset V_t \subset \C^{t+1} \]
with $\dim V_i =i$ such that
\[ D(\Pe(V_j),X_i) \geq - c_{10} \deg X_i, \quad i=0,\ldots,t,j=1, \ldots,t. \]
Consequently, by Proposition \ref{linausw},
\[ D(\Pe(V_j),X) \geq - (c_{10}+C) \deg X, \quad j=1,\ldots,t. \]
\end{Satz}

\proof
For $i$ arbitrary, by Proposition \ref{Dfunk1},
\[ D(X_i, \Pe(V_i)) = D(V_{X_i},V_i), \]
and by the previous Proposition, there is a flag such that this is greater or 
equal $- c_{10}(i) \deg X_i$.

\vspace{2mm}

{\sc Proof of Theorem \ref{Punktraum}}
By the Corollary, there are spaces $\Pe(V_{p-2})\subset\Pe(V_{p-1})$
of dimensions $p-2,p-1$ respectively such that with $X_\theta$ as in (\ref{XW}),
\[ D(X_\theta,\Pe(V_{p-2})) \geq - c_{10}\deg X_\theta = -c_{10}\deg X, \quad
   D(X,\Pe(V_{p-1})) \geq - c_{10}\deg X. \]
Choosing $\Pe(F)$
as any subspace containing $\Pe(V_{p-1})$ as well as $\theta$,
the claim follows from Proposition \ref{linausw}.

\vspace{2mm}

2. Let $\Pe(F)$ be admissible for $X$ and $\theta$.
By Proposition \ref{raumraum}.2,
\[ D(\theta,X) = D^{\Pe(F)}(\theta,X.\Pe(F))  + D(X,\Pe(F))-D(X_\theta,\Pe(F)) 
   + c_{11} \deg X. \]
The claim thus follows from Remark \ref{remark} and the fact that
$\Pe(F)$ is admissible.

\section{Joins}

The proof of Theorem \ref{bezout} will be using the construction of
the join of cycles $\CX,\CY \subset \Pe(E) = \Pe^t$ in projective space, 
which is defined as follows.
Let $\pi_i: \Pe(E) \to \Pe(E) \times \Pe(E), i=1,2$ be the canonical 
embeddings, $\pi:\Pe(E\oplus E) \to \spec(\Z)$
the structure maps, and
$p_i: \Pe(E) \times \Pe(E) \to \Pe(E), i=1,2$ the canonical projections. Then
$F := p_1^* O(-1) \oplus p_2^* O(-1)$ is a subbundle of
\[ \pi^* (E\oplus E) = p_1^* \pi^* E \oplus p_2^* \pi^* E \]
over $\Pe(E) \times \Pe(E)$.
The inclusion defines a map from the projective bundle \\
$\Pe_{\Pe(E) \times \Pe(E)}(F)$ 
to $\Pe_{\Pe(E) \times \Pe(E)}(\pi^*(E \oplus E)) = 
\Pe(E \oplus E) \times \Pe(E) \times \Pe(E)$, hence a map $g$ to
$\Pe^{2t+1} = \Pe(E\oplus E)$; the bundle map 
$f:\Pe_{\Pe(E) \times \Pe(E)}(F) \to \Pe(E) \times \Pe(E)$ is flat. Hence,
for $\CX,\CY \in \Pe(E)$, the expression
\begin{equation} \label{joindef}
\CX\#\CY := g_*f^*(\CX \times \CY) 
\end{equation}
is well defined and is called the join of $\CX$ and $\CY$. We have

%Alternatively, let $\CX$ is given by the homogenous polynomials
%$s_1,\ldots s_n$ in the variables $x_0, \ldots x_t$. These can be viewed
%%as polynomials $s'_1, \ldots s'_n$ in the variables
%$x_0, \ldots x_t,y_0,\ldots y_t$ for $\Pe(E\oplus E)$.
%Likewise one has homogeneous polynomials
%$t_1', \ldots t_m'$ for $\CY$. Then $\CX\#\CY$ is the
%subvariety of $\Pe(E\oplus E)$ defined by the homogenous polynomials
%$s'_1, \ldots s_n', t_1', \ldots t_m'$.

\satz{Proposition} \label{join}
The degree and height of the join compute as
\[ \deg(X\#Y) = \deg(X) \deg(Y), \quad \mbox{and} \]
\[ h(\CX\#\CY) = \deg X h(\CY) + \deg Y h(\CX). \]
\end{Satz}

\proof
\cite{BGS}, Proposition 4.2.2.

\satz{Proposition} \label{joinD}
With $p+r \geq t+1, q+s \geq t+1$ let $X,Y,Z,W$ be effective cycles in
$\Pe^t_\C$ of codimensions $p,q,r,s$ respectively such that the $X$ and $Z$ as
well as $Y$ and $W$ intersect properly. Then, $X\#Y$ and $Z\#W$ intersect
properly, and the algebraic distance $D(X\#Y,Z\#W)$ computes as
\[ D(X \# Y, Z \# W) = \deg X \deg Z D(Y,W) + \deg Y \deg W D(X,Z)+ \]
\[ (\sigma_{2t+1}-2\sigma_t) \deg X \deg Y \deg Z \deg W. \]

Let $\theta \in \Pe^t(\C)$ be a point neither contained in the support
of $X$ nor in that of $Y$, and $X_\theta,Y_\theta$ the varieties as in 
(\ref{XW}). Further, $(\theta,\theta) \in \Pe^{2t+1}$ the intersection
of $\theta \# \theta$ with the diagonal, and $(X\# Y)_{\theta,\theta}$
as in (\ref{XW}). Then,
\[ D((X \# Y)_{(\theta,\theta)},Z \# W) = \deg X \deg Z D(Y_\theta,W)+
   \deg Y \deg W D(X_\theta,Z) + \]
\[ (\sigma_{2t+1}-2\sigma_t) \deg X \deg Y \deg Z \deg W. \]

\end{Satz}

\proof
It suffices to prove the Proposition for cycles that are irreducible varieties.
Assume first that $X,Y,Z,W$ are obtained by base extension from
subvarieties $\CX,\CY,\CZ,\CW$ of $\Pe^t_\Z$.
By Remark \ref{remark},
\[ D(X \# Y, Z \# W) = h((\CX \# \CY) . (\CZ \# \CW)) - \]
\[ \deg (X \# Y) h(\CZ \# \CW) - \deg (Z \# W) h(\CX \# \CY) +
     \sigma_{2t+1} \deg (X \# Y) \deg (Z \# W), \]
\[  h((\CX . \CZ) \# (\CY \# \CW)) -
     \deg (X \# Y) h(\CZ \# \CW) - \deg (Z \# W) h(\CX \# \CY) + \]
\[   \sigma_{2t+1} \deg (X \# Y) \deg (Z \# W), \]
which by the previous Proposition and Remark \ref{remark} equals
\[ \deg X \deg Z (\deg Y h(\CW) + \deg W h(\CY) +
                  D(Y,W) - \sigma_t \deg Y \deg Z) + \]
\[ \deg Y \deg W (\deg X h(\CZ) + \deg Z h(\CX) + D(X,Z) - 
                 \sigma_t \deg X \deg Z) - \]
\[ \deg X \deg Y (\deg Z h(\CW) + \deg W h(\CZ)) -
   \deg Z \deg W (\deg X h(\CY) + \deg Y h(\CX)) + \]
\[ \sigma_{2t+1} \deg X \deg Y \deg Z \deg W = \]
\[ \deg X \deg Z D(Y,W) + \deg Y \deg W D(X,Z)+
   (\sigma_{2t+1}-2\sigma_t) \deg X \deg Y \deg Z \deg W. \]
Since the the subvarieties $\CX \subset \Pe^t(\Z)$ form a dense subset
of the subvarieties $X \subset \Pe^t(\C)$, and the algebraic distance is
a continuous function, the equation holds for arbitrary $X,Y,Z,W$.

Since $(X\# Y)_{\theta,\theta} = g_* f^*(X_\theta \times Y_\theta)$, 
by Lemma \ref{Dfunk1},
\[ D((X \# Y)_{(\theta,\theta)},Z \# W) = D(X_\theta \# Y_\theta,Z \# W), \]
which by the above equals
\[ D((X \# Y)_{(\theta,\theta)},Z \# W) = \deg X \deg Z D(Y_\theta,W)+
   \deg Y \deg W D(X_\theta,Z) + \]
\[ (\sigma_{2t+1}-2\sigma_t) \deg X \deg Y \deg Z \deg W. \]

\vspace{2mm}

Let $X,Y$ be cycles of pure codimension $p,q$ with $p+q \geq t+1$. 
They do intersect
properly, do not intersect respectively, if
if and only if $X \# Y$ does intersect $\Pe(\Delta)$ properly, does
not intersect $\Pe(\Delta)$.
Thus, one may define

\satz{Definition}
For $p + q \leq t+1$ define the algebraic distance
\[ \bar{D}(X,Y) := D(\Pe(\Delta),X\#Y)), \]
For $p+1\geq t+1$ define 
\[ \bar{D}_G(X,Y) := D_G(\Pe(\Delta),X\#Y), \quad
   \bar{D}_\infty(X,Y) := D_\infty(\Pe(\Delta),X\#Y), \]
\[ \bar{D}_{Ch}(X,Y) := D_{Ch}(\Pe(\Delta),X\#Y). \]
\end{Satz}

As $(X + X') \# Y = (X \# Y) + (X' \# Y)$ it immediately follows form
Lemma \ref{addit}, that for $p+q \geq t+1$ the maps
\[ \bar{D}_0, \bar{D}_{Ch}: Z^p_{eff}(\Pe^t) \times Z^q_{eff}(\Pe^t) \to \R \]
are bilinear.

\satz{Proposition} \label{neu}
There exist constants $c,c_0,c_\infty,c_{Ch}$ only depending on $p,q$, 
and $t$ such that in the situation of the Definition,
\[ D(X,Y) = \bar{D}(X,Y) + c \deg X \deg Y, \]
and for $Y = \Pe(W)$ a projective subspace,
\[ \bar{D}_0(X,\Pe(W)) = D_0(\Pe(W),X)+ c_0 \deg X, \]
\[ \bar{D}_\infty(X,\Pe(W)) = D_\infty(\Pe(W),X) + c_\infty \deg X, \]
\[ \bar{D}_{Ch}(X,\Pe(W)) = D_{Ch}(\Pe(W),X) + c_{Ch} \deg X. \]
That is, the above Definition of algebraic distances coincide with
the old definition modulo constants times $\deg X$.
\end{Satz}

\proof
The proof will be given in a different paper (\cite{liouville}).

\vspace{2mm}

Let $\CX,\CY$ be irreducible subschemes of $\Pe^t$ whose generic fibre is
not empty, and $[x],[y]$ closed points of $X_\C,Y_\C$ represented by vectors
$x,y \in\C^{t+1}$ of length one. Then the closed points of the join
$x \# y \subset (\CX \# \CY)_\C$ are the points
\[ g(f^{-1}([x],[y])) = g (\lambda x,\mu y) = [(\lambda x,\mu y)], \]
with $\lambda,\mu \in \C$.

\satz{Lemma} \label{joindist}
Let $x,y,\theta$ be vectors of length one in $\C^{t+1}$ with 
$[x] \neq [\theta] \neq [y]$. Then,
\[ \mbox{min} (|[\theta] ,[x]|,|[\theta],[y]|) \leq 
   |[x] \# [y], ([\theta],[\theta])| \leq
   \mbox{max} (|[\theta] ,[x]|,|[\theta],[y]|). \]
\end{Satz}

\proof
Since
$|\la([\theta],[\theta])\ra,[x]\#[y]| = |\la([\theta],[\theta])\ra,\la [v]\ra|$
where $[v]$ is the point in $x\#y$ with minimal distance to 
$[(\theta,\theta)]$. Hence, it suffices to prove
\[ \mbox{min} (|[\theta] ,[x]|,|[\theta],[y]|) \leq |[v], (\theta,\theta)| \leq
   \mbox{max} (|[\theta] ,[x]|,|[\theta],[y]|). \]
As any point $[w]$ in $x \# y \in \Pe^t(\C)$ may be written as 
$[\lambda (x,0) + \mu (0,y)]$ with $\lambda, \mu \in \C$, and
\[ |(\lambda x,\mu y),(\theta,\theta)|^2 = 
   1 - \frac{|\la (\lambda x, \mu y) |(\theta,\theta)\ra |^2}
   {2(|\lambda|^2+|\mu|^2)}; \]
further, 
$|x,\theta|^2 = 1-|\la x|\theta\ra|^2, |y,\theta|^2= 1-|\la y|\theta\ra|^2$,
we have to show that
\[ \mbox{min}(|\la x|\theta\ra|^2,|\la y|\theta\ra|^2) \leq
   {\mbox{sup} \atop \lambda,\mu} 
   \frac{|\la \lambda x, \mu y |(\theta,\theta)\ra |^2}
   {2(|\lambda|^2+|\mu|^2)} \leq
   \mbox{max}(|\la x|\theta\ra|^2,|\la y|\theta\ra|^2). \]
Firstly,
\[ \frac{|\la (\lambda x, \mu y) |(\theta,\theta)\ra |^2}
   {2(|\lambda|^2+|\mu|^2)} \leq
   \frac{|\lambda^2 \la x|\theta \ra^2| + 
         2 |\lambda  \mu  \la x|\theta \ra \la y|\theta \ra|^2 +
         |\mu^2 \la y|\theta \ra^2|}
   {2(|\lambda|^2+|\mu|^2)} \leq \]
\[  \frac{|\lambda|^2 + 2|\lambda \mu| + |\mu|^2}
   {2(|\lambda|^2+|\mu|^2)} \;
   \mbox{max}(|\la x|\theta\ra|^2,|\la y|\theta\ra|^2) \leq
   \mbox{max}(|\la x|\theta\ra|^2,|\la y|\theta\ra|^2), \]
as $|\lambda|^2 + 2|\lambda \mu| + |\mu|^2 \leq
2(|\lambda|^2+|\mu|^2)$, whence the second inequality.

For the first inequality, choose
\[ \lambda_0 = 
   \sqrt{\frac{|\la x|\theta\ra|^2}{|\la x|\theta\ra|^2 +
   |\la y|\theta\ra|^2}}, \quad
   \mu_0 =
   \sqrt{\frac{|\la y|\theta\ra|^2}
   {|\la x|\theta\ra|^2 +|\la y|\theta\ra|^2}}. \]
Then,
\[  \frac{|\la (\lambda_0 x, \mu_0 y) |(\theta,\theta)\ra |^2}
   {2(|\lambda_0|^2+|\mu_0|^2)} =
   \frac12 |\la(\lambda_0 x,\mu_0 y),(\theta,\theta)\ra|^2 =
   \frac{(|\la x|\theta \ra|^2+|\la y|\theta \ra|^2)^2}
            {2(|\la x|\theta \ra|^2+|\la y|\theta \ra|^2)} = \]
\[ \frac{|\la x|\theta \ra|^2+|\la y|\theta \ra|^2}{2} \geq 
   \mbox{min}(|\la x|\theta \ra|^2,\la y|\theta \ra|^2), \]
whence the first inequality.

%We will need a slight extension of Lemma \ref{raumraum}
%
%\satz{Lemma} \label{raumraumj}
%Let $X,Y$ be effective cycles of pure codimension $p$ and $q$ in 
%$\Pe^t$, and $\theta$ a point not contained in either of their
%supports. With $\Delta \subset \C^{t+2}$ the diagonal,
%and $V \subset \C^{t+2}$ such that $\theta \# \theta = \Pe(V)$,
%let $(\Pe(\Delta + V), \pi, \psi_\lambda,(X \# Y)_\lambda)$ be as in
%(\ref{deform}), and $\Pe(W)$ any subspace of $\Pe(\Delta)$ of
%codimension $r \leq t+1-p-q$. Then, for any $\lambda_1, \lambda_2 \in \C$,
%\[ D(\Pe(W),(X \# Y)_{\lambda_1}) - D(\Pe(W),(X \# Y)_{\lambda_2}) = \]
%\[ D(\Pe(\Delta + V),(X \# Y)_{\lambda_1}) -
%   D(\Pe(\Delta+V),(X \# Y)_{\lambda_2})= \]
%_\[ D(\Pe(\Delta),(X \# Y)_{\lambda_1})-D(\Pe(\Delta),(X \# Y)_{\lambda_2}). \]
%\end{Satz}

%\proof
%Looking at the proof of Lemma \ref{raumraum}, it only has to be shown
%that $(X \# Y)_\lambda . \Pe(\Delta + V) = (X \# Y) . \Pe(\Delta + V)$
%for all $\lambda \in \C$, and this is easily checked. 

\section{Proofs for the main results}

By Theorem \ref{egal},  the results are true
for one of the algebraic distances $D_G, D_\infty, D_{Ch}$ if it holds for any
of the others, only with different constants $c,c',e,e'$. 

Remember, that for a subvariety $\CX \subset \Pe^t_\Z$, the height
$h(\CX)$ is defined for $\CX$ over $\Z$, the degree is defined for the
base extension $X = \CX_\Q$, and the algebraic distance to some
point $\theta \in \Pe^t(\C)$ is defined for the $\C$- valued points of $\CX$,
denoted $X_\infty$, or $X$ if clear from the context.

\proof {\sc of Proposition \ref{haupt1}}
Let $p$ be the codimension of $\CX$, and
$\Pe(F)$ a $c_{10}$-admissible subspace of dimension $p$ for $X$ and $\theta$.
Then,
\[ \deg X \log |\theta,X| \leq \sum_{x \in X\cap \Pe(F)} n_x \log |\theta,x|=
   D_{\Pe(F)}(\theta,X), \]
which by Theorem \ref{Punktraum} is less or equal 
$D(\theta,X) + (c_{10} + e_1) \deg X$ proving the first inequality.

For the second inequality, let
$x_0 \in \mbox{supp}(X(\C))$ be such that $|\Pe(W),X|= |\Pe(W),x_0|$.
By Corollary \ref{adm}, there is a subspace
$\Pe(V) \subset \Pe^t$ of dimension $p-2$ such that
$D(\Pe(V),X_\theta) \geq - c_{10} \deg X$.
With $\Pe(F)$ a space of dimension $p$ containing $\Pe(V)$ as well as $\theta$
and $x_0$ 
and intersecting $X$ as well as $X_\theta$ properly, Proposition \ref{linausw}
implies $D(\Pe(F),X_\theta) \geq - \left(c_{10}+C\right)\deg X$.
Further by Proposition \ref{raumraum}.2,
\[ D(\theta,X) = D^{\Pe(F)}(\theta,X.\Pe(F)) + D(\Pe(F),X) - D(\Pe(F),X_\theta)
   + c_{11} \deg X, \]
which by the above and Remark \ref{remark} is less or equal
\[ D^{\Pe(F)} (\theta,X.\Pe(F)) + 
   \left(\bar{c}_1+ \frac{c_{10}}{d(p-2,t)} + C + c_{11}\right) \deg X, \]
which in turn by Proposition \ref{abst} equals
\[ \sum_{x \in X\cap \Pe(F)} n_x \log |\theta,x| +
   \left(\bar{c}_1+ \frac{c_{10}}{d(p-2,t)} + C + c_{11}+c_2+c_4\right)\deg X = \]
\[ n_{x_0} \log |\theta,x_0| + \sum_{(x \in X)\setminus \Pe(F)} n_x 
   \log |\theta,x|+
   \left(\bar{c}_1+ \frac{c_{10}}{d(p-2,t)} + C + c_{11}+c_2+c_4\right)
   \deg X \]
\[ \leq \log |\theta,X| +
   \left(\bar{c}_1+ \frac{c_{10}}{d(p-2,t)} + C + c_{11}+c_2+c_4\right)\deg X, \]
finishing the proof.

\proof {\sc of Proposition \ref{haupt2}}
Remember that in case of codimension one $D_\infty =D_G=D_{Ch}$.
To deduce the first inequality, firstly, by \ref{hsch}.\@1, and (\ref{hpr}),
\[ h(\di\; f) = D \sigma_t + \int_{\Pe^t(\C)} \log|f| \mu^t. \]
Further, by Lemma \ref{normrelpr},
\[  \int_{\Pe^t(\C)} \log|f| \mu^t \leq
    \frac12\log \int_{\Pe^t(\C)} |f|^2 \mu^t = \log |f|_{L^2} \]
The two formulas together imply the first formula.

For the second formula, 
\[ D(\theta,\di \;f)  = \log |f_{\theta}| - \int_{\Pe^t_\C} \log |f| \mu^t, \]
which by Proposition \ref{hsch} and (\ref{hpr}) equals
\[ \log |(f|_\theta)| - h(\di \;f) + D \sigma_t. \]
%As by (\ref{Levint}) and the definition of the Stoll numbers,
%\[ \int_{\Pe^t(\C)} \Lambda_\theta \mu = \sum_{n=1}^t \frac{1}{n} =
%   (\sigma_t - \sigma_{t-1}), \]
%we get
%\[ -\frac{1}{2} \int_{\di(f_D)(\C)} \Lambda_\theta + h(\di(f_D)) =
%    \log |\la f_D |\theta \ra| + \frac D2 (\sigma_t + \sigma_{t-1}). \] 

\proof {\sc of Theorem \ref{bezout}}
With $c=c_{10} +C$,where $c_{10}$ and $C$ are the constants from
Theorem \ref{Punktraum} and Proposition \ref{linausw}, 
let $\Pe(F),\Pe(F')$ be $c$-admissible subspaces of dimensions
$p,q$ for  $X$ and $\theta$, $Y$ and $\theta$ respectively, and let
$x_1,\ldots,x_{\deg X},\\ y_1,\ldots,y_{\deg Y} \in \Pe^t(\C)$ such that
\[ |\theta,x_1| \leq |\theta,x_2| \leq \cdots \leq |\theta,x_{\deg X}|, \quad
   |\theta,y_1| \leq \cdots \leq |\theta,y_{\deg Y}|, \]
and
\[ X . \Pe(F) = \sum_{i=1}^{\deg X} x_i, \quad 
   Y . \Pe(F')= \sum_{i=1}^{\deg Y} y_i. \]
By Theorem \ref{Punktraum},
\begin{equation} \label{absa}
\sum_{i=1}^{\deg X} \log |\theta,x_i| \leq
D(\theta, X) + (c+e_2) \deg X,
\end{equation}
and
\begin{equation} \label{absb}
\sum_{i=1}^{\deg Y} \log |\theta,y_i| \leq
D(\theta, Y) + (c+e_2) \deg Y.
\end{equation}

Define the function $f_{X,Y}$ subject to requirement that for every
$T \in \underline{\deg X + \deg Y}$ and $(\nu,\kappa) = f_{X,Y}(t)$
the conditions
\[ \nu + \kappa = t, \quad \mbox{and} \quad
   |\theta,x_\nu| \leq |\theta,y_{\kappa+1}|, \quad 
   |\theta,y_{\kappa}| \leq |\theta,x_{\nu + 1}| \]
hold. 

The proof will be given in two steps. 

\vspace{2mm}

1.
\[ \nu \kappa \log |\theta,X+Y| + D_G((\theta,\theta),X\#Y) + 
   h(\CX\#\CY) \leq \]
\[ \kappa D_G(\theta,X) + \nu D_G(\theta,Y) 
   +\deg Y h(\CX) +\deg X h(\CY) + \]
\[ (2\sigma_t-\sigma_{2t+1}+3c+e_1 + e_2) \deg X \deg Y, \]

\vspace{2mm}

2.
\[ D_\infty(\theta,X . Y) + h(\CX . \CY) \leq \]
\[ D_\infty((\theta,\theta),X \# Y) + h(\CX \# \CY) + 
   \left(\frac{3t+2-p-q}2 \log 2+\sigma_t-\sigma_{2t+1}\right) \deg X \deg Y. \]

\vspace{2mm}

As, by Theorem \ref{Punktraum},
$D_\infty((\theta,\theta),X \#Y) \leq 
D_G((\theta,\theta),X \# Y) + c_3 \deg X \deg Y$,
and $D_G(\theta,X.Y) \leq D_\infty(\theta,X.Y)$,
the two inequalities together imply the claim with
$e = 2\sigma_t-2\sigma_{2t+1}+3c+e_1 + e_2 + \frac{3t+2-p-q}2 \log 2$.

\vspace{2mm}

1. Write $D(\cdot,\cdot)$ for $D_G(\cdot,\cdot)$.
By Proposition \ref{joinD},
\[ D(X \# Y, \Pe(F) \# \Pe(F')) = 
   \deg X D(Y,\Pe(F')) + \deg Y D(X,\Pe(F)) + \]
\[ (\sigma_{2t+1}-2\sigma_t) \deg X \deg Y \geq 
   (\sigma_{2t+1}-2\sigma_z - 2c) \deg X \deg Y, \]
and
\[ D(X \# Y, \Pe(F) \# \Pe(F')) \leq 
   \bar{c}_1 \deg X \deg Y. \]
Further,
\[ D((X \# Y)_{(\theta,\theta)},\Pe(F) \# \Pe(F')) = \]
\[ \deg X D(Y_\theta,\Pe(F')) + \deg Y D(X_\theta,\Pe(F)) +
   (\sigma_{2t+1}-2\sigma_t) \deg X \deg Y \geq \]
\[ (\sigma_{2t+1}-2\sigma_t-2c) \deg X \deg Y, \]
and 
\[  D((X \# Y)_{(\theta,\theta)},\Pe(F) \# \Pe(F')) \leq \bar{c}_1 \deg X \deg Y. \]

Hence, $\Pe(F)\#\Pe(F)'$ is a $(2\sigma_t-\sigma_{2t+1}+2c)$-admissible subspace
for $X\#Y$ and $(\theta,\theta)$. Since
\[ (\Pe(F) \# \Pe(F')) . (X \# Y) = 
   \sum_{i=1}^{\deg X} \sum_{j=1}^{\deg Y}
   x_i \# y_j, \]
Theorem \ref{Punktraum} and Proposition \ref{abst}.2 imply
\begin{eqnarray} 
D((\theta,\theta), X \# Y) 
               & \leq & 
    D^{\Pe(F) \# \Pe(F')}((\theta,\theta),(X\#Y).(\Pe(F) \# \Pe(F')))\nonumber\\ &+&
             (2\sigma_t-\sigma_{2t+1} +2c+e_1) \deg X\nonumber \\
             &=&
          \sum_{x,y} n_y n_y \log |(\theta,\theta),x\#y|  \nonumber
          \\ \label{haupteins}&+&
                 (2\sigma_t-\sigma_{2t+1} +2c+e_1) \deg X\#Y. 
\end{eqnarray}
Next, since the logarithm of the Fubini-Study metric is nonpositive, 
\[ \nu \kappa \log|\theta,X+Y| +
   \sum_{i=1}^{\deg X} \sum_{j=1}^{\deg Y} 
   \log |x_i\# y_j,(\theta,\theta)| \leq \]
\[ \nu \kappa \log |\theta,X + Y| + 
\sum_{i=1}^\nu \sum_{j=1}^\kappa \log |x_i \#y_j ,(\theta,\theta)| + \]
\begin{equation} \label{abk} 
\sum_{i=1}^\nu \sum_{j=\kappa+1}^{\deg Y} \log |x_i \# y_j,(\theta,\theta)| +
\sum_{i=\nu+1}^{\deg X} \sum_{j=1}^\kappa \log |x_i \# y_j,(\theta,\theta)|. 
\end{equation}
Now, since 
$\log |\theta,X + Y| \leq \mbox{min}(\log|\theta,x|,\log|\theta,y|)$,
for any $x \in \mbox{supp} \;(X), y \in \mbox{supp} \; (Y)$,
and by Lemma \ref{joindist} also
$\log |(\theta,\theta), X \#Y| \leq \mbox{max}(\log|\theta,x|,\log|\theta,y|)$,
for all $i \leq \nu,j \leq \kappa$, the inequality
\[ \log|\theta,X+Y| + \log |x_i \# y_j,(\theta,\theta)| \leq
   \log \mbox{min}(|x_i,\theta|,|y_j,\theta|) + \]
\[ \log \mbox{max}(|x_i,\theta|,|y_j,\theta|) = 
   \log |x_i,\theta| + \log |y_j,\theta| \]
holds.
Hence, the sum of the first two summands on the right hand side 
of (\ref{abk}) is less or equal than
\begin{equation} \label{abd}
\kappa \sum_{i=1}^\nu \log|x_i,\theta|+ \nu \sum_{j=1}^\kappa \log |y_j,\theta|.
\end{equation}
For $i \leq \nu$ and $j \geq \kappa+1$ we have 
$|x_i,\theta| \leq |y_j,\theta|$, consequently, by Lemma \ref{joindist}
$|x_i \# y_j,(\theta,\theta)| \leq |\theta,y_j|$.
Using this, and the analogous inequality for
$i \geq \nu+1,j \leq \kappa$, the sum of the third and fourth summand of 
(\ref{abk}) is less or equal
\begin{equation} \label{abf}
\sum_{i=1}^\nu \sum_{j=\kappa+1}^{\deg Y}  \log |y_j,\theta| +
   \sum_{i=\nu+1}^{\deg X} \sum_{j=1}^\kappa \log |x_i,\theta| =
\nu \sum_{j=\kappa+1}^{\deg Y} \log |y_j,\theta| + 
\kappa \sum_{i=\nu+1}^{\deg X} \log |x_i,\theta|. 
\end{equation}
Hence, (\ref{abk}) is less or equal than the sum of (\ref{abd}) 
and (\ref{abf}), which equals
\[ \kappa \sum_{i=1}^{\deg X} \log |x_i,\theta| +
   \nu \sum_{j=1}^{\deg Y} \log |y_j,\theta|, \]
which in turn, by (\ref{absa}), and (\ref{absb}) is less or equal 
\[ \kappa D(\theta,X) + \nu D(\theta,Y) + 
   (e_2+c)(\nu \deg Y + \kappa \deg X). \]
Together with (\ref{haupteins}), this gives
\[ \nu \kappa \log |\theta,X+Y| + D((\theta,\theta), X \# Y) \leq 
   \kappa D(\theta,X) + \nu D(\theta,Y) + \]
\[ ((2\sigma_t-\sigma_{2t+1}+2c+e_1) \deg X \deg Y + 
   (e_2+c)(\kappa \deg Y + \nu \deg X) \leq \]
\[ \nu D(\theta,Y) + \kappa D(\theta,X) + 
   (2\sigma_t-\sigma_{2t+1}+3c+e_1+e_2) \deg X \deg Y, \]
%which, by Proposition \ref{sum} is less or equal
%\[ \nu (D_0(\theta,X) - D(\Pe(F),X)) + 
%   \kappa (D_0(\theta,Y)- D_0(\Pe(F'),Y)) + (2e_3-e_2) \deg X \deg Y. \]
%This in turn, by Proposition \ref{absch} is greater or equal
%\[ \nu D_0(\theta,X) + \kappa D_0(\theta,Y) +
%   (\nu \deg Y+ \kappa \deg X)(c_4+c_3+c_3-c_1) + 
%   (2e_3-e_2) \deg X \deg Y \leq \]
%\[ \nu D_0(\theta,X) + \kappa D_0(\theta,Y) + (2(c_4+c_3+c_3-c_1) +2e_3-e_2)
%   \deg X \deg Y. \]
Adding the equation $h(\CX \# \CY) = \deg(\CX) h(\CY) + \deg(\CY) h(\CX)$
of proposition  \ref{join}.\@2 to this inequality leads the desired inequality.

\vspace{2mm}

2. Let $\Delta \subset \C^{2t+2}$ be the diagonal, and $i:\Pe^t \to \Pe^{2t+1}$ 
the inclusion
\[ i: \Pe^t \to \Pe(\Delta) \subset \Pe^{2t+1}, \quad [v] \mapsto [(v,v)]. \]
Then, $(X \# Y) . \Pe^{\Delta} = i(X.Y)$, and since the restriction
of $\overline{O_{\Pe^{2t+1}}(1)}$ to $\Pe^t=\Pe(\Delta)$ equals
$\overline{O_{\Pe^t}(1)}$ with the norm multiplied with $\sqrt{2}$, 
we have $h((\CX\# \CY).\Pe(\Delta)) = h(\CX.\CY)+ \frac{2t+2-p-q}2\log 2$.

Let $\Pe(V) \subset \Pe^t$ be a subspace of dimension
$p+q-1$ containing $\theta$ such that $D(\Pe(V),X . Y)$ is maximal, 
i.e.\@ $D_\infty(\theta,X.Y) = D(\Pe(V),X.Y)$.

Since $i:\Pe^t \to \Pe(\Delta)$ is an isometry, we get
\[ D_\infty(\theta,X.Y) = D(\Pe(V),X.Y) = 
   D^{\Pe(\Delta)}(i(\Pe(V)),(X\#Y).\Pe(\Delta)), \]
which by Proposition \ref{raumraum}.1 equals
\[ D(i(\Pe(V)),X\#Y) - D(\Pe(\Delta),X\#Y), \]
which in turn by Remark \ref{remark} equals 
\[ D(i(\Pe(V)),X\#Y) + h(\CX \# \CY) + \deg(X\#Y) h(\Pe(\Delta)) -
   h(\Pe(\Delta).\CX\#\CY) - \]
\[ \sigma_{2t+1} \deg X\#Y. \]
Since  $h(\Pe(\Delta)) = \frac t2 \log 2+ \sigma_t$, by (\ref{hpr}),
and $h((\CX\# \CY).\Pe(\Delta)) = h(\CX.\CY)+\frac{2t+2-p-q}2\log 2$ from above,
\[ D_\infty(\theta,X.Y) + h(\CX . \CY) = \]
\[ D(i(\Pe(V)),X\#Y) + h(\CX\#\CY) + 
   \left(\frac{3t+2-p-q}2\log 2+\sigma_t-\sigma_{2t+1}\right) \deg X \deg Y, \]
and the trivial estimate $D(i(\Pe(V)),X\#Y) \leq D_\infty((\theta,\theta),X\#Y)$
implies
\[ D_\infty(\theta,X.Y) +h(\CX.\CY) \leq \]
\[ D_\infty((\theta,\theta),X\#Y) +h(\CX\#\CY) + 
   \left(\frac{3t+2-p-q}2\log 2-\sigma_{2t+1}\right) \deg X \deg Y, \]
as was to be proved

\proof {\sc of Corollary \ref{haupt3}}
1. The way $f_{X,Y}$ is defined in the proof of Theorem  \ref{bezout} above,
it follows that $(\nu,\kappa)$ in Theorem \ref{bezout} can be chosen such that 
$\nu = 1$. Then,
\[ \kappa \log |\theta,X+Y| + D(\theta,X.Y) + 
   h(\CX . \CY) \leq \]
\[ \kappa D(\theta,X) + D(\theta,Y) +
   \deg Y h(X) + \deg X  h(Y) + e \deg X \deg Y. \]
As by Theorem \ref{bezout}, 
\[ \kappa D(\theta,X) \leq \kappa \log |\theta, X| + \kappa c(p,t) \deg X, \]
in either of the two cases
$|\theta,X| \leq |\theta,Y|$ or $D(\theta,X) \leq \log |\theta,Y|$ one has
\[ \kappa D(\theta,X) \leq \kappa \log |\theta,X+Y| + \kappa c(p,t) \deg X. \]
Inserting this into the first inequality and subtracting 
$\kappa \log |\theta,X+Y|$ from the resulting formula
gives the claim with $e' = c+e$.

\vspace{2mm}

2. Again, in Theorem \ref{bezout} $(\nu,\kappa)$ can be either chosen equal
to $(0,1)$ or equal to $(1,0)$. Without loss of generality, assume
$(\nu,\kappa) = (0,1)$. Then,
\[ D(\theta,X.Y) + h(\CX.\CY) \leq D(\theta,X) + 
   \deg Y h(X) + \deg X  h(Y) + e \deg X \deg Y \leq \]
\[ \max(D(\theta,X),D(\theta,Y)) +\deg Y h(X) + \deg X  h(Y) + 
   e \deg X \deg Y. \]

%\vspace{2mm}

%{\bf Remark:}
%The algebraic distance $D_\bullet(\Pe(F),X)$ 
%has been essentially divided into two 
%Summands in Proposition \ref{sum}, and essentially reduced to one
%of these summands in Theorem \ref{Punktraum}. The proof of
%the metric B\'ezout Theorem presented here estimates only this essential
%summand in the algebraic distance of the joint to the same summand for the
%two cycles. It is actually possible to do the same thing for the 
%second summand. This would give the same metric B\'ezout Theorem,
%but with a somewhat better constant $d$.

\satz{Proposition} For $X,Y$ arbitrary effective cycles in $\Pe^t(\C)$, and 
$\theta \in \Pe^t(\C)$ a point not contained in $supp X \cup supp Y$
there is a function 
$f_{X,Y} :\underline{\deg X+\deg Y} \to \underline{\deg X} \times
\underline{\deg Y}$ with the same properties as in Theorem \ref{bezout}
such that for every $T \in f_{X,Y} :\underline{\deg X+\deg Y}$ and
$(\nu,\kappa) =f_{X,Y}(T)$ the inequality
\[ \nu \kappa \log |\theta,X+Y| + D(\theta, X . Y)+ D(X,Y) \leq
   \kappa D(\theta,X) + \nu D(\theta,Y) + \bar{e} \deg X \deg Y. \]
Further, if $|\theta,X+Y|= |\theta,X|$, or $D(\theta,X) \leq \log |\theta,Y|$,
then
\[ D(\theta, X.Y) + D(X,Y) \leq
   D(\theta,Y) + \bar{e}' \deg X \deg Y,  \]
and in general
\[ D(\theta,X.Y) + D(X,Y) \leq \max(D(\theta,X),D(\theta,Y)) +
   \bar{e}'\deg X\deg Y. \]
If $X$, and $Y$ have pure complementary dimension, then $D(\theta,X.Y) = 0$,
and the above implies the logarithmic triangle inequality
\[ D(X,Y) \leq \mbox{max}(D(\theta,X), D(\theta,Y)) + 
   \bar{e}' \deg X \deg Y \]
holds.
\end{Satz}

\proof
Just repeat the proofs of Theorem \ref{bezout} and Corollary
\ref{haupt3} without using the fact
\[ D(\Pe(\Delta),\CX \# \CY) = \]
\[ h((\CX \# \CY) . \Pe(\Delta)) - h(\CX \# \CY) - 
   \deg (X \#Y) h(\Pe(\Delta)) + \sigma_{2t+1} \deg (X\#Y), \]
and instead use that $D(\Pe(\Delta),\CX \# \CY)$, and $D(X,Y)$ only
differ by a constant times $\deg X \deg Y$ by
Proposition \ref{neu}.

%The first part of the proof used the triangle inequality for the
%Fubini-Study metric on $\Pe^t_\C$, and the fact that the algebraic distance
%$D$ is expressed in terms of the Fubini-Study metric. The second part
%uses the fact that the algebraic distance $D$ is closely related to an
%algebraic distance $c_3$ described by
%Green forms. The idea to proove the
%metric B\'ezout Theorem without the factors $2t$ and $2$, would be to find 
%a definition of algebraic distance that combines both of these properties.


\begin{thebibliography}{0mm}

\bibitem[BGS]{BGS} Bost, Gillet, Soul\'e: Heights of projective varieties and
positive Green forms. JAMS 7,4 (1994)

\bibitem[Fu]{Fu} W.\@ Fulton: Intersection Theory. 2nd edition, Springer 1998.

\bibitem[GS1]{GS1} H.\@ Gillet, C.\@ Soul\'e: Arithmetic intersection theory.
Publications Math.\@ IHES 72 (1987) 243-278

\bibitem[GS2]{GS2} H.\@ Gillet, C.\@ Soul\'e: Characteristic classes for 
algebraic vector bundles with hermitian metric I,II. 
Annals of Mathematics 131 (1990), 163-238.

\bibitem[LR]{LR} M.\@ Laurent, D.\@ Roy: Criteria of algebraic independence
with multiplicities and approximation by hypersurfaces.

\bibitem[Mai]{Mai} V.\@ Maillot: Un Calcul de Schubert arithm\'etique.
Duke Math.\@ J.\@ 80,1 (1995), 195-221

\bibitem[Ma1]{App2} H.\@ Massold: Diophantine Approximation on varieties II:
  Explicit estimates for arithmetic Hilbert functions. arxiv: 0711.1667

\bibitem[Ma2]{App3} H.\@ Massold: Diophantine Approximation on varieties III:
  Approximation of non-algebraic points by algebraic points. arxiv: 0711.3645

\bibitem[Ma3]{App4} H.\@ Massold: Diophantine Approximation on varieties IV:
  Derivated algebraic distance and derivative metric B\'ezout Theorem.
  arxiv: 0901.3889

\bibitem[Ma4]{App5} H.\@ Massold: Diophantine Approximation on varieties V:
  Algebraic independence criteria.
  arxiv: 1001.1534

%\bibitem[Ma5]{mahler} H.\@ Massold: Mahlerklassifiezierung von Punkten auf
%  algebraischen Variet\"aten.

\bibitem[Ma6]{liouville} H.\@ Massold: Liouville's Theorem for subvarieties
 of arbitrary dimension. To appear.

%\bibitem[Ma7]{trans}H.\@ Massold: Estimates for Transcendence degrees of fields
%   generated by values of exponentials of algebraic groups I: Commutative 
%   groups. To appear.


%\bibitem[Ma4]{Ma4} H.\@ Massold: Transcendence properties of values of
%   E-functions on transcendental points. To appear.

\bibitem[Nes]{Nes} Y.\@ V.\@ Nesterenko: On the measure of algebraic
independence of the values of Ramanujan functions. Trudy Matematicheskovo
Instituta imany V.\@ A.\@ Steklov 218 (1997), 299-334. English translation in:
Proj.\@ Steklov Inst.\@ Math.\@ 218 (1997), 294-331.

\bibitem[Ph1]{Ph} P.\@ Philippon: Sur des hauteurs alternatives I,
  Math.\@ Ann.\@ 289, (1991) 255-283; II, Ann.\@ Inst.\@ Fourier 44 (1994)
1043-1065; III, J.\@ Mathe.\@ Pures Appl.\@ 74 (1995) 345-365.

\bibitem[Ph2]{Ph1} P.\@ Philippon: Approximations alg\'ebriques des points
dans les espaces projectifs I. J.\@ Number Theory 81, No.\@ 2 234-253 (2000)

\bibitem[RW]{RW} D.\@ Roy, M.\@ Waldschmidt: Approximation diophantienne
et ind\'ependance alg\'ebrique de logarithmes. Ann.\@ sc.\@ de l'ENS
30,6 (1997)

\bibitem[SABK]{SABK} C.\@ Soul\'e, D.\@ Abramovich, J.\@-F.\@ Burnol,
               J.\@ Kramer: Lectures on Arakelov Geometry. 
               Cambridge University Press 1992 

\end{thebibliography}
\end{document}